\DocumentMetadata{
    pdfversion=1.7,
    pdfstandard=a-3b,
}

\documentclass[conference,10pt]{IEEEtran}

\IEEEoverridecommandlockouts

\usepackage[margin=1.5cm]{geometry}

\usepackage{mathtools}

\usepackage{xcolor}

\usepackage{bm}

\usepackage{siunitx}

\usepackage{ifluatex}

\usepackage{subcaption}

\usepackage{graphicx}

\ifluatex
    \usepackage{fontspec}
    \usepackage{unicode-math}
    \NewDocumentCommand\vectorStyle{m}{\symbf{#1}}
    \NewDocumentCommand\matrixStyle{m}{\symbf{#1}}
\else
    \ifxetex
        \usepackage{fontspec}
        \usepackage{unicode-math}
        \NewDocumentCommand\vectorStyle{m}{\symbf{#1}}
        \NewDocumentCommand\matrixStyle{m}{\symbf{#1}}
    \else
        \NewDocumentCommand\vectorStyle{m}{\bm{#1}}
        \NewDocumentCommand\matrixStyle{m}{\bm{#1}}
    \fi
\fi

\usepackage{amssymb}

\usepackage[sorting=none,style=ieee]{biblatex}

\usepackage{balance}

\NewDocumentCommand\position{}{\vectorStyle{x}}

\NewDocumentCommand\positions{}{\matrixStyle{X}}

\NewDocumentCommand\currentDensity{}{\vectorStyle{j}}

\NewDocumentCommand\primaryCurrentDensity{}{{\currentDensity_{\mathrm{p}}}}

\NewDocumentCommand\electricPotential{}{u}

\NewDocumentCommand\transferMatrix{}{\matrixStyle{T}}

\NewDocumentCommand\leadFieldMatrix{}{\matrixStyle{L}}

\DeclarePairedDelimiter\args{(}{)}

\NewDocumentCommand\Hdiv{}{H(\mathrm{div})}

\DeclarePairedDelimiter{\norm}{\lVert}{\rVert}

\NewDocumentCommand\Nof{m}{N_{\text{#1}}}

\NewDocumentCommand\nodeN{}{\Nof{n}}

\NewDocumentCommand\electrodeN{}{\Nof{e}}

\NewDocumentCommand\sourceN{}{\Nof{s}}

\NewDocumentCommand\numberSetStyle{m}{\mathbb{#1}}

\NewDocumentCommand\realNumbers{}{\numberSetStyle{R}}

\NewDocumentCommand\functionsFromTo{mm}{{#1}^{#2}}

\NewDocumentCommand\meshNodePosition{}{\position_{\text{n}}}

\NewDocumentCommand\electrodePosition{}{\position_{\text{e}}}

\NewDocumentCommand\sourcePosition{}{{\position_{\text{s}}}}

\NewDocumentCommand\sourcePositions{}{\positions_{\text{s}}}

\DeclareMathOperator\opEMD{EMD}

\NewDocumentCommand\cost{}{\matrixStyle{C}}

\NewDocumentCommand\flow{}{\matrixStyle{F}}

\NewDocumentCommand\estimate{m}{{#1}^{\operatorname{est}}}

\begin{document}

\title{Forward--Inverse Interplay in FEM-Based EEG Source Imaging: Distributional Signatures of Advanced Source Models and Inverse Solvers%
\thanks{This work was supported by the Research Council of Finland through the Flagship of Advanced Mathematics for Sensing, Imaging and Modelling (FAME), 2024--2031.}
}

\author{
\IEEEauthorblockN{Santtu Söderholm, Joonas Lahtinen, and Sampsa Pursiainen}
\IEEEauthorblockA{\textit{Mathematics Research Center, Faculty of Information Technology and Communication Sciences} \\
\textit{Tampere University}\\
Korkeakoulunkatu 1, 33720 Tampere, Finland}
}

\maketitle

\begin{abstract}
Electroencephalography (EEG) source imaging aims to infer brain activity from electrical potentials generated by current sources measured on the scalp. This is a fundamentally ill-posed problem because there are multiple orders of magnitude more variables related to the source space than there are sensors. The result depends strongly on two things: the forward model and the inverse method. In this work, we study how these two parts work together. We focus not only on where the activity is located, but also on how the reconstructed activity is distributed in space. We suggest that different source models create different signatures in the reconstructed activity. We use realistic head models and compute forward solutions with the finite element method using Zeffiro Interface and DUNEuro. We test different source models, including 2 implementations of a divergence-conforming model, and one implementation of Local subtraction approach. For inverse methods, we use advanced methods such as standardized hierarchical adaptive L1 regression (SHAL1R), standardized Kalman filtering (SKF), and classical Dipole Scanning (DS). To understand the complex interplay between the forward and inverse approaches, we analyze the inverse source localization results using distributional quantitative measures, including Earth Mover's Distance and depth bias scatter plot, and qualitatively assess the amplitude distribution and focality. The results show that there is a strong dependence between the choice of source model and the success rate of a given inverse method: a source model that corresponds well with a single point-like source is a good match with an inverse method that presupposes such a source.
\end{abstract}

\begin{IEEEkeywords}
electroencephalography, finite element method, inverse problems, source imaging
\end{IEEEkeywords}

\section{Introduction}

Electroencephalography (EEG) source imaging aims to estimate brain activity from signals measured on the scalp. This is an ill-posed inverse problem because the dimensionality of the source current space is practically always orders of magnitude larger than that of the measurement space, leading to non-uniqueness in an inverse solution. An inverse solution is also sensitive to noise and modelling errors in the forward phase. In practice, the result depends on two main components: the lead field matrix $\leadFieldMatrix\in\functionsFromTo\realNumbers{\electrodeN\times\sourceN}$, which describes how brain activity from fixed $\sourceN$ source locations generates signals at $\electrodeN$ number of sensors, and the inverse method, which reconstructs the sources from these measurements. A lot of progress has been made in developing inversion methods tailored to source imaging. New variations of the standardized method, sLORETA \cite{pascual2002sloreta}, have been introduced in recent years \cite{SSLOFO2005Liu,lahtinen2024shalpr,lahtinen2024standardized}. The appeal of standardization stems from its theoretical ability to reconstruct sources without depth bias from noiseless data \cite{Pascual2007sqrtm,lahtinen-standardization-2024}.

Recent studies have also reaffirmed that modelling choices in the forward problem can significantly affect the reconstruction outcome. In particular, preprocessing steps such as peeling \cite{soderholm2024effects}, which restrict the positions where synthetic sources can be placed in a discretized computational domain where conductivity jumps exist, can have a notable effect on numerical EEG forward solution \cite{wolters-etal-2007} and therefore the corresponding source reconstruction results. This further emphasizes that forward modelling choices influence not only accuracy but also the spatial characteristics of the reconstructed activity. Finite-Element-Method (FEM)-based forward modelling provides a flexible framework for improving modelling accuracy by incorporating realistic head geometry and tissue conductivity distributions \cite{he2020zeffiro,schrader2021duneuro}. Within this framework, the choice of source model significantly impacts the reconstruction. For example, divergence-conforming $\Hdiv$ models (a quadratic extension of the linear Whitney basis), partial integration, and St.~Venant approaches have been shown to exhibit different trade-offs between focality and numerical accuracy \cite{miinalainen2019realistic}. Importantly, these differences are reflected not only in localization error but also in the spatial distribution of the reconstructed activity.

Recent developments, such as the Local subtraction approach, provide a mathematically sound and computationally efficient way to model singular dipolar sources in FEM \cite{holtershinken2025local}. Combining those with advanced inverse methods enables reconstructing weakly distinguishable activity \cite{lahtinen2024shalpr,lahtinen2024standardized}, e.g., deep subcortical sources, while also motivating a comprehensive analysis of EEG source imaging dynamics that goes beyond pointwise localization accuracy.

In this paper, we utilize the unbiasedness of standardized estimation and dipolar fitting in Dipole Scan (DS) to assess differences arising from the $\Hdiv$ source model implementation of Zeffiro Interface and two alternative implementations provided by DUNEuro: Whitney source model with face-intersecting and edgewise basis functions and Local subtraction. The results show that disagreement with the source model and assumptions of the inversion methods can cause deviations from the ideal results.

\section{Background}\label{sec.background}

\begin{figure}
    \centering
    \includegraphics[width=0.3\linewidth]{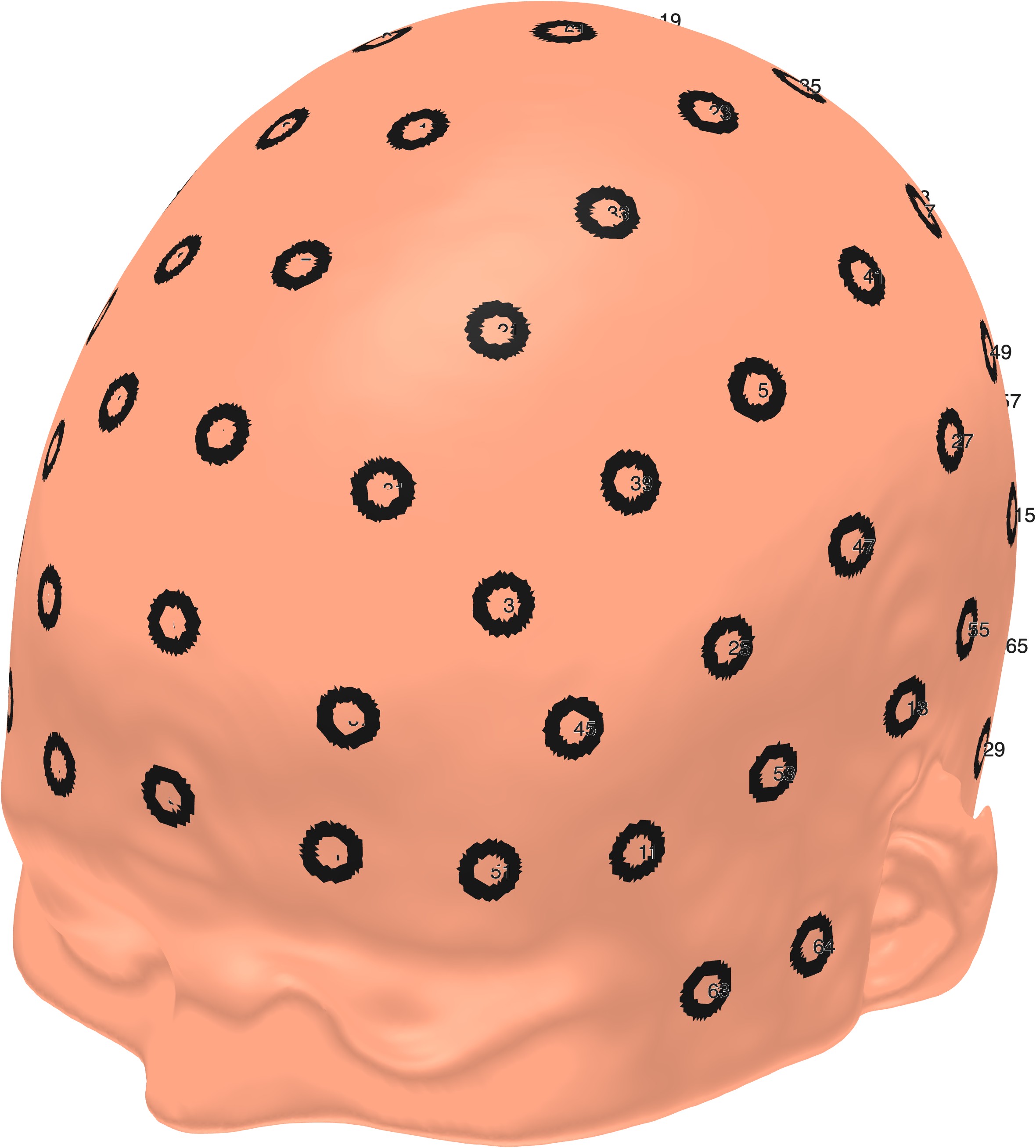}
    \caption{\small Electrodes positioned on the scalp of the ICBM152 2009a head model \cite{fonov-2009}, according to the international 10--10 electrode positioning standard \cite{Knosche-Haueisen-2022}.}
    \label{fig.scalp.electrodes}
\end{figure}

Attempting to determine brain activity from voltage measurements is a \num{100}-year-old endeavor, originating in Germany \cite{vergani-2024}. Originally, the means of performing the measurements were rather invasive, requiring the use of depth electrodes to make the signals discernible with the technology of the day. While depth electrodes are still in use today when accurate measurements are required \cite{missey-etal-2026}, typical modern EEG measurements are performed transcranially, with the electrodes placed on the scalp \cite{Knosche-Haueisen-2022}, with an example seen in Figure~\ref{fig.scalp.electrodes}.

Mathematical modelling of the EEG sources also has a history of decades \cite{clark-and-plonsey-1966}. With forward modelling, a rough progression has been from the use of simple spherical models \cite{ary-etal-1981, cuffin2001spherical} towards realistic domains acquired via segmentation of MRI images \cite{cuffin2001realistic}, with boundary element methods (BEM) \cite{demunck-1992,schlitt-1995} and more recently Finite Element Methods (FEM) \cite{wolters-etal-2007} working as the computational approaches. Inverse solvers originally relied on a form of dipole fitting, where a semi-analytical forward solution is compared with a measured one and the dipole model parameters were optimized such that the difference between the forward solution and the measured signal was minimized \cite{ary-etal-1981}. More recent inverse solvers either rely on the idea of regularizing a lead field matrix $\leadFieldMatrix$, or use Bayesian modelling with more complicated a priori assumptions and iterative algorithms in the seeking of maximum a posteriori that maximizes the log-posterior density function \cite{kaipio2006statistical,Knosche-Haueisen-2022}.

\section{Methods}\label{sec.methods}

Lead field matrices are maps from sources to EEG electrodes  \cite{Knosche-Haueisen-2022}. They are computed using the FEM-based transfer matrix approach implemented in DUNEuro \cite{schrader2021duneuro} and Zeffiro Interface \cite{he2020zeffiro}. Inverse reconstructions of brain activity are produced with solvers found in Zeffiro Interface. The freely available and realistic population-based ICBM152 2009a multicompartment head model \cite{fonov-2009} is used to account for anatomical structure and conductivity variations. Our model of the electrodes is the point electrode model (PEM) \cite{hanke-2011,agsten2018electrodes}.

When constructing a lead field, we consider different \emph{source models}, which correspond to means of interpolating a transfer matrix $\transferMatrix\in\functionsFromTo\realNumbers{\nodeN\times\electrodeN}$ from the nodes of a finite element mesh $\meshNodePosition$ to given source positions $\sourcePosition$ \cite{wolters2004transfer,Knosche-Haueisen-2022}. Chosen source models include the face-intersecting and edge-wise Whitney basis functions, the divergence-conforming $\Hdiv$ formulation \cite{miinalainen2019realistic}, and the Local subtraction approach \cite{holtershinken2025local}.

Here, we consider the differences, advantages, and disadvantages of the source models from an inversion perspective. Namely, we focus on the smoothness, spread, and depth bias of inversion estimations with various methods. In our experiment I, we estimate a cortical source under inversion crime to assess how smooth and spread the estimation we get, and compare the difference of the estimation distribution visually to the estimations obtained with Local subtraction. In this experiment, we use an anisotropic model of the brain conductivity. In experiment II, we compute source estimates at different depths from the inner skull layer of the upper scalp (in the vicinity of the sensor surface). Sources were positioned by randomly placing \num{100} Cartesian sources per \qty{5}{\milli\meter} interval on the S- or $z$-axis of the RAS coordinate system, in relation to the lowest $z$-coordinate of the set of electrodes. The positions were in the range of \num 0--\qty{60}{\milli\meter} relative height, resulting in a total of \num{1200} sources. In this experiment, we use an anisotropic lead field to generate the synthetic data and an isotropic model in inversion to avoid inversion crime while keeping the same mesh used in head model creation.

The source models are compared using the depth bias scatter plot introduced in \cite{elvetun2025depthbias} and the \emph{Earth Mover's Distance} (EMD) to measure the spread of the reconstructions. Earth Mover's Distance is a measure of similarity between two distributions, and is computed between our sets of original source positions $\sourcePositions$, and their reconstructed estimates $\estimate\sourcePositions$
\cite{rubner-etal-1998-emd}:
\begin{equation}\label{eq.emd}
    \opEMD
    \args * {
        \sourcePositions,
        \estimate\sourcePositions
    }
    =
    \frac{
        \sum_i
        \sum_j
        \cost_{i,j}
        \flow_{i,j}
    }{
        \sum_i
        \sum_j
        \flow_{i,j}
    }
    \,.
\end{equation}
Here $i$ corresponds to a true source position $\sourcePosition_i\in\sourcePositions$ and $j$ to a reconstructed source position $\estimate\sourcePosition_j\in\estimate\sourcePositions$.
The quantity $\cost_{i,j}$ is the cost of moving source $i$ to $j$, or the distance between them, and $\flow_{i,j}\leq\norm{\primaryCurrentDensity_i}, \norm{\estimate\primaryCurrentDensity_j}$ is a non-negative flow restricted by the true and estimated primary current densities $\primaryCurrentDensity_i$ and $\estimate\primaryCurrentDensity_j$, that minimizes the expression in the numerator of \eqref{eq.emd}.

For source reconstruction, we use standardized distributed methods including sLORETA \cite{pascual2002sloreta}, the sparsity-promoting SHAL1R \cite{lahtinen2024shalpr}, and spatiotemporal SKF \cite{lahtinen2024standardized} with classical Dipole Scanning (DS) \cite{fuchs-1998,neugebauer-etal-2022}, from which well-established sLORETA and DS serve as benchmark methods. Because DS is observed to be sensitive mainly to the difference between forward and inverse models \cite{Lahtinen2023}, we omit it in Experiment II. SKF is omitted in Experiment II due to its dependence on the evolution model's agreement with the dynamics of time-dependent source activity.

The motivation for displaying the depth-bias scatter plot for standardized methods stems from the theoretical result that the estimated source location should be unbiased, meaning we obtain perfect agreement between the true source depth and the estimated depth \cite{Pascual2007sqrtm}. We assume that greater deviation from the perfect agreement indicates elevated ambiguity among sources, which is undesirable. In EMD terms, a lower value indicates a smaller spread. The evolution of EMD, along with the depth of the true source, is also an important indicator.

\section{Results}\label{sec.results}

\begin{figure}
    \centering
    \def\figWidth{0.2\linewidth}
    \def\sepWidth{0.1\linewidth}
    \def\colorbarHeight{2cm}
    \begin{subfigure}[b]{\figWidth}
        \includegraphics[width=\linewidth]{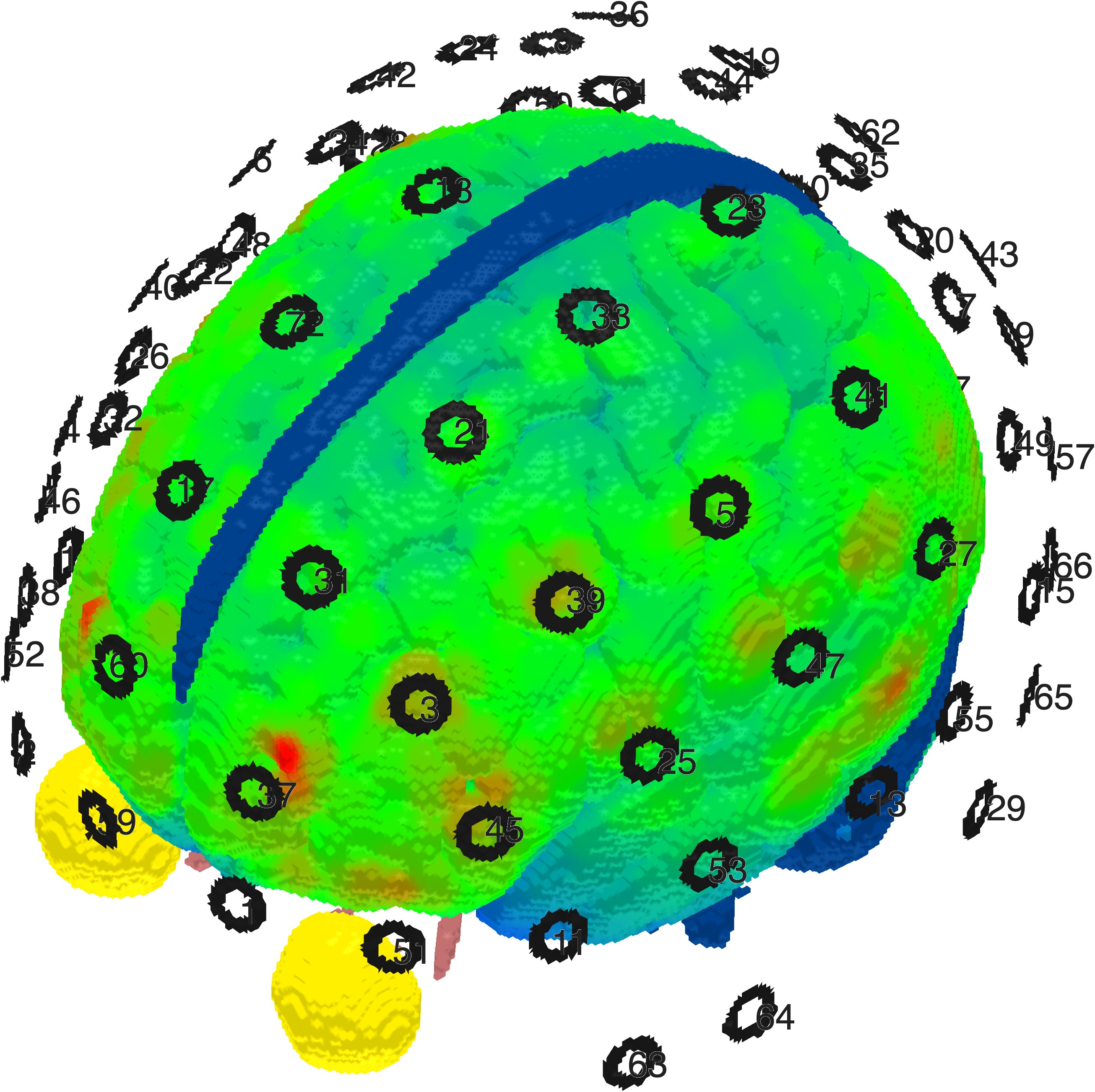}
        \caption{}
        \label{fig.duneuro.L.vecnorm.superficial}
    \end{subfigure}
    \hspace\sepWidth
    \begin{subfigure}[b]{\figWidth}
        \includegraphics[width=\linewidth]{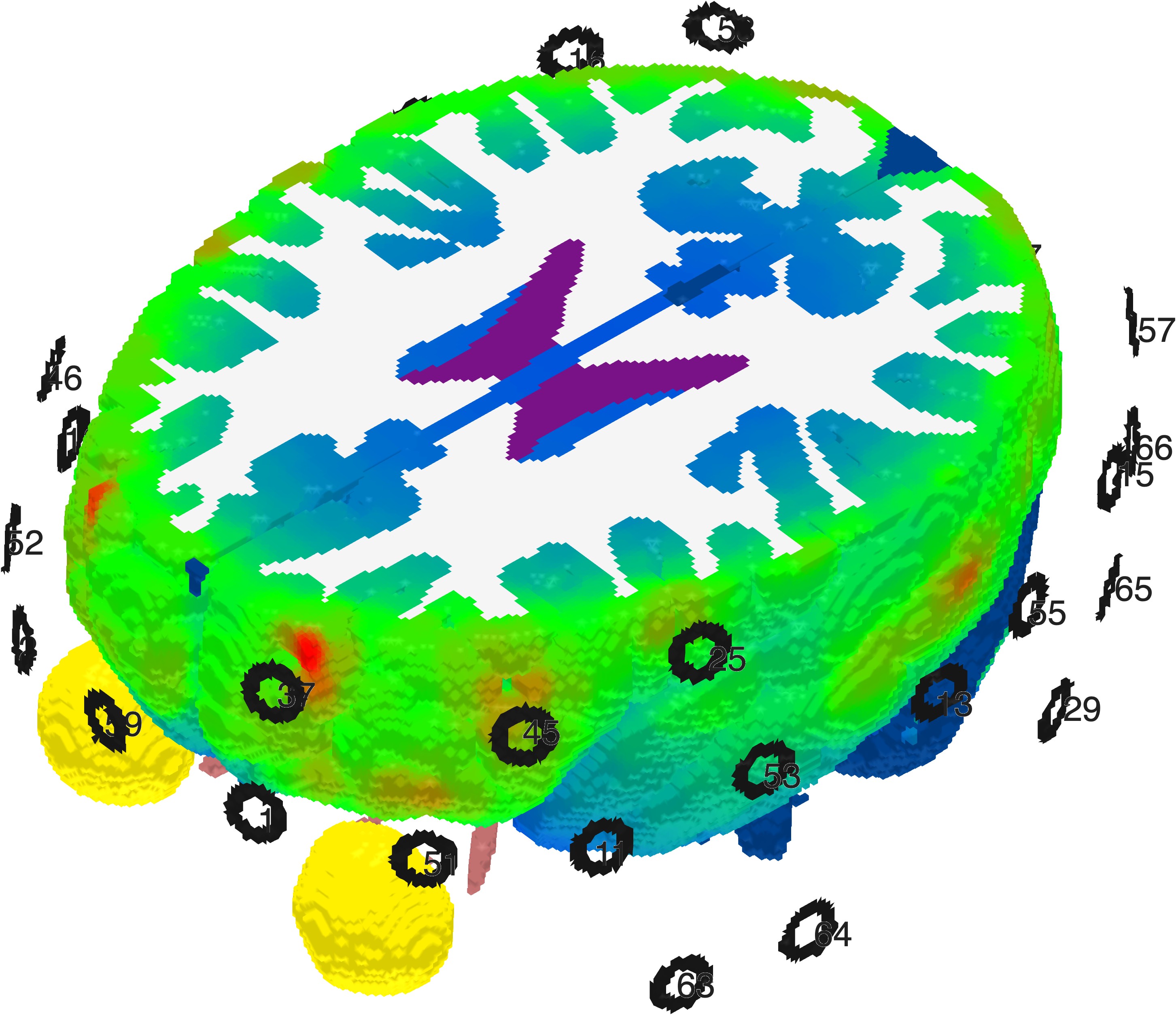}
        \caption{}
        \label{fig.duneuro.L.vecnorm.deep}
    \end{subfigure}
    \hspace\sepWidth
    \begin{minipage}[b]{0.1\linewidth}
        \includegraphics[height=\colorbarHeight]{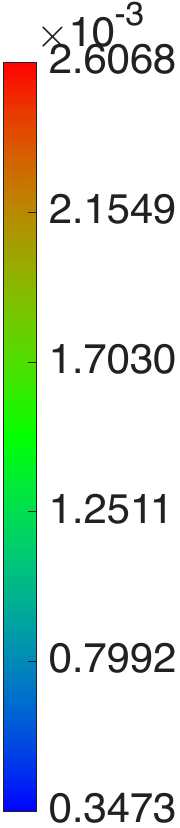}
    \end{minipage}

    \vspace*{0.1cm}

    \begin{subfigure}[b]{\figWidth}
        \includegraphics[width=\linewidth]{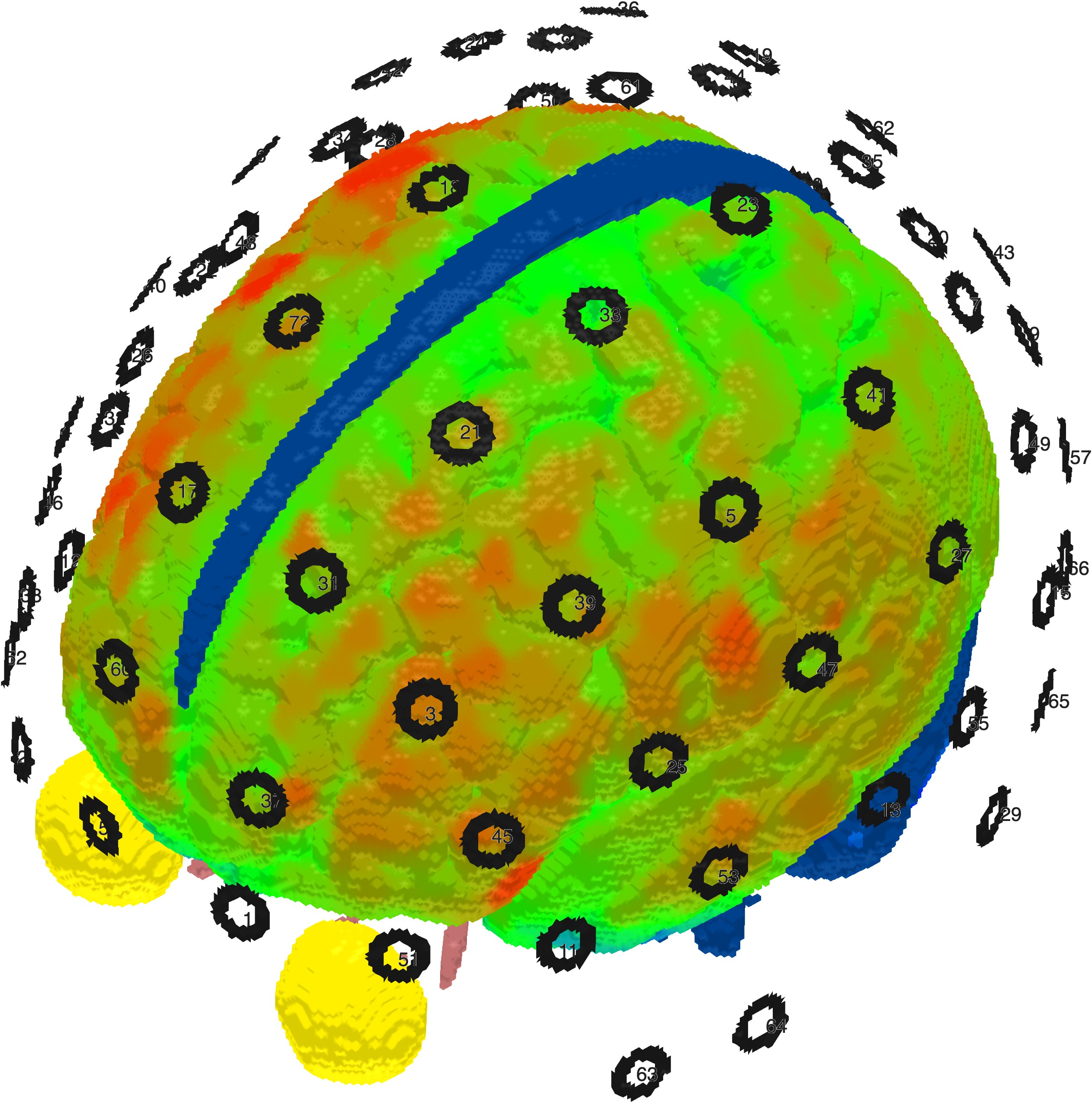}
        \caption{}
        \label{fig.zeffiro.L.vecnorm.superficial}
    \end{subfigure}
    \hspace\sepWidth
    \begin{subfigure}[b]{\figWidth}
        \includegraphics[width=\linewidth]{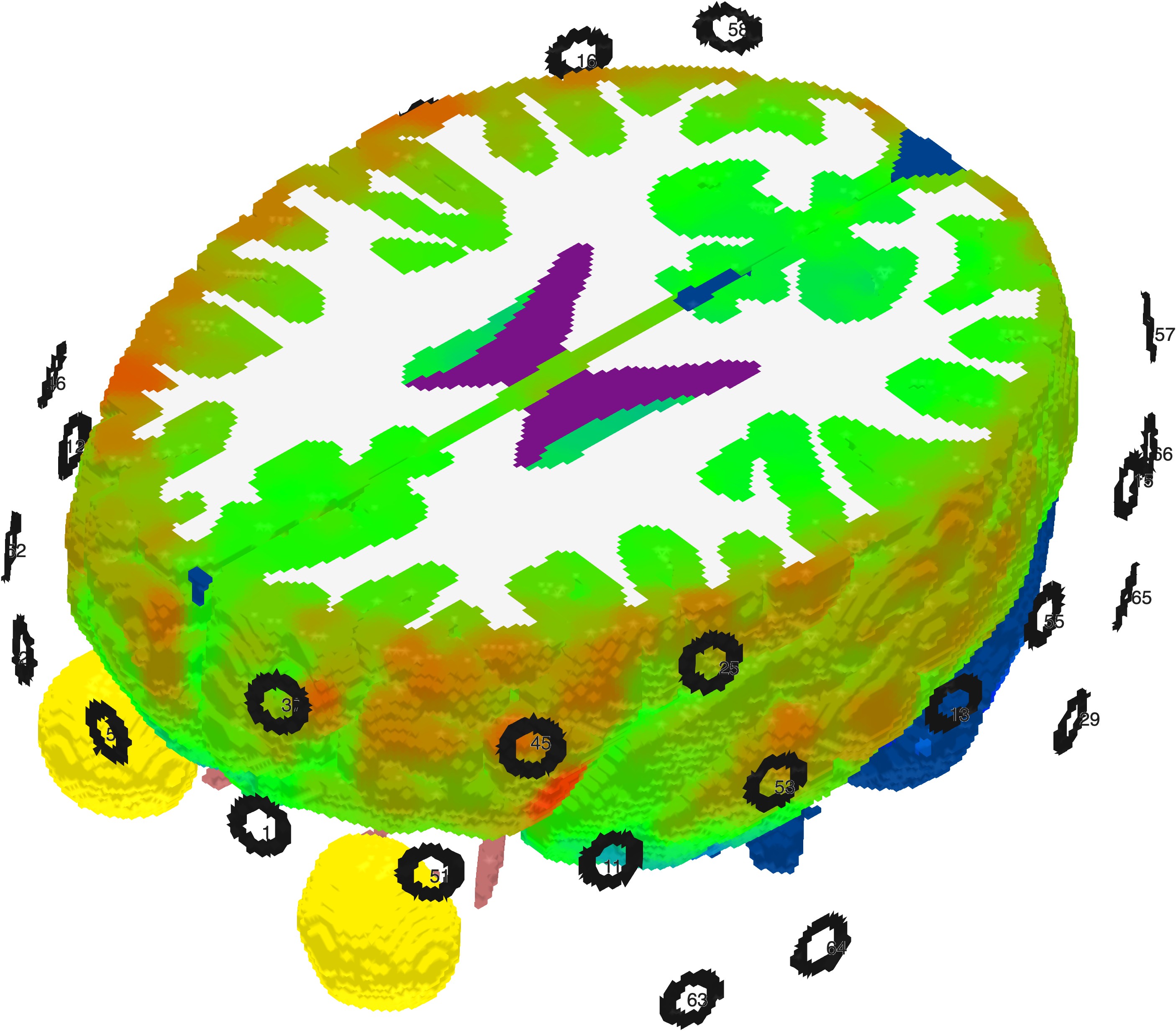}
        \caption{}
        \label{fig.zeffiro.L.vecnorm.deep}
    \end{subfigure}
    \hspace\sepWidth
    \begin{minipage}[b]{0.1\linewidth}
        \includegraphics[height=\colorbarHeight]{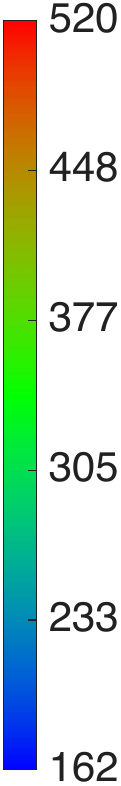}
    \end{minipage}

    \caption{\small A display of the column norms $\norm\leadFieldMatrix_2$ of lead fields $\leadFieldMatrix$ produced with Whitney-based DUNEuro (\ref{fig.duneuro.L.vecnorm.superficial}, \ref{fig.duneuro.L.vecnorm.deep}) and $\Hdiv$-based Zeffiro Interface (\ref{fig.zeffiro.L.vecnorm.superficial}, \ref{fig.zeffiro.L.vecnorm.deep}). Subfigures~\ref{fig.duneuro.L.vecnorm.superficial} and \ref{fig.zeffiro.L.vecnorm.superficial} display a cortical view while \ref{fig.duneuro.L.vecnorm.deep} and \ref{fig.zeffiro.L.vecnorm.deep} show a deeper cross-section in the RA-plane of the Right-Anterior-Superior (RAS) coordinate system, roughly at the height of the thalamus. The DUNEuro field is strong near the EEG electrodes depicted as black circles, and decays smoothly when moving away from them. The Zeffiro field displays higher irregularity and had to be low-pass filtered up to the \qty{95}{\percent} quantile to get rid of an unwarranted field maximum occurring behind the yellow eyeballs. The DUNEuro field is also orders of magnitude smaller in value.}
    \label{fig.L.vecnorms}
\end{figure}

Figure~\ref{fig.L.vecnorms} (a--d) displays examples of column norms of a lead field computed with DUNEuro and Zeffiro Interface. The DUNEuro lead fields have been computed using the Whitney basis functions, while Zeffiro uses $\Hdiv$ interpolation in the source space. Interpolation of the vector basis functions from the mesh nodes and element faces and edges to the synthetic source positions is done using Position-Based Optimization (PBO) \cite{pursiainen-etal-2016}, which allows the source positions to deviate from the mesh structure itself, instead of being anchored to the mesh geometry and introducing mesh-based bias to the solution. Separate source generation also allows one to limit the size of the lead field used in the inverse phase, making computations more feasible.

The norm of the Whitney-based $\leadFieldMatrix$ behaves as expected in Figure~\ref{fig.L.vecnorms} (a--b): electric potential $\electricPotential$ falls off with the inverse distance from a point source. The fact that the field described by $\transferMatrix$ is constructed by first placing sources at the PEM electrode positions $\electrodePosition$, then solving for a potential field from them to the mesh nodes $\meshNodePosition$ and finally using Helmholtz reciprocity \cite{vanrumste2004reciprocity,gross2023reciprocity} or the transposition of $\transferMatrix$ to produce a mapping from $\meshNodePosition$ to $\electrodePosition$. This is then interpolated to the source positions $\sourcePosition$ to produce $\leadFieldMatrix$.
The $\Hdiv$-based field in Figure~\ref{fig.L.vecnorms} (c--d) is more irregular, and the decay of the lead field column norm is weaker from the vicinity of the electrodes towards the thalamus compared to Whitney.

\begin{figure}
    \def\imageWidth{0.18\linewidth}
    \def\inverseMethodNameWidth{0.08\linewidth}
    \def\colorbarWidth{0.06\linewidth}
    \tiny
    \centering
    \begin{minipage}[b]{0.8\linewidth}
            \begin{minipage}[b]{\inverseMethodNameWidth}
            \rotatebox{90}{\hspace{0.05cm} sLORETA}
        \end{minipage}\begin{minipage}[b]{\imageWidth}
        \begin{center}
            Whitney
        \end{center}
            \includegraphics[trim={6cm 2cm 6.5cm 2cm},clip,width=1\linewidth]{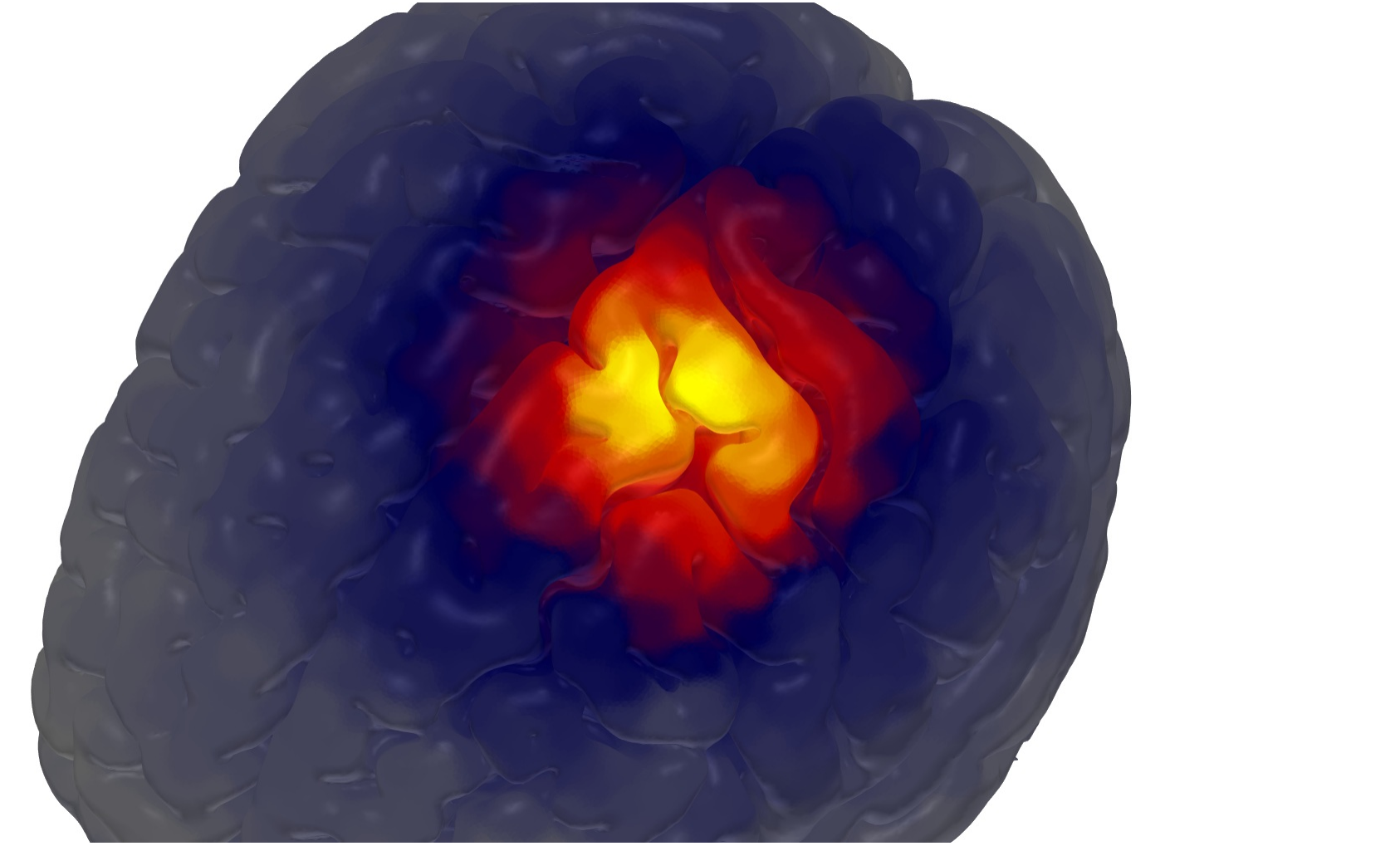}
        \end{minipage}\begin{minipage}[b]{\imageWidth}
        \begin{center}
            Local subtraction
        \end{center}
            \includegraphics[trim={6cm 2cm 6.5cm 2cm},clip,width=1\linewidth]{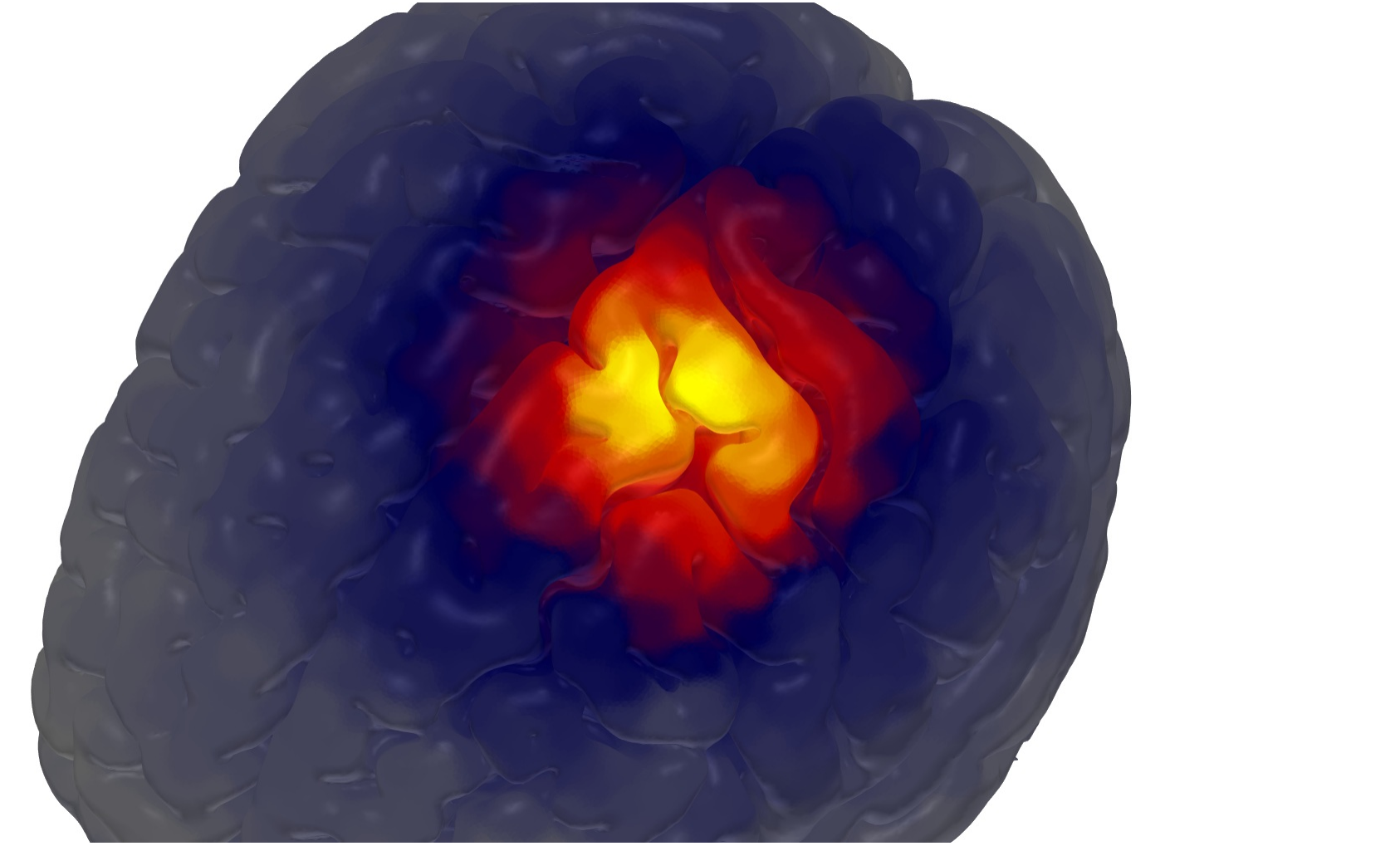}
        \end{minipage}%
        \begin{minipage}[b]{\imageWidth}
        \begin{center}
            $\Hdiv$
        \end{center}
            \includegraphics[trim={6cm 2cm 6.5cm 2cm},clip,width=1\linewidth]{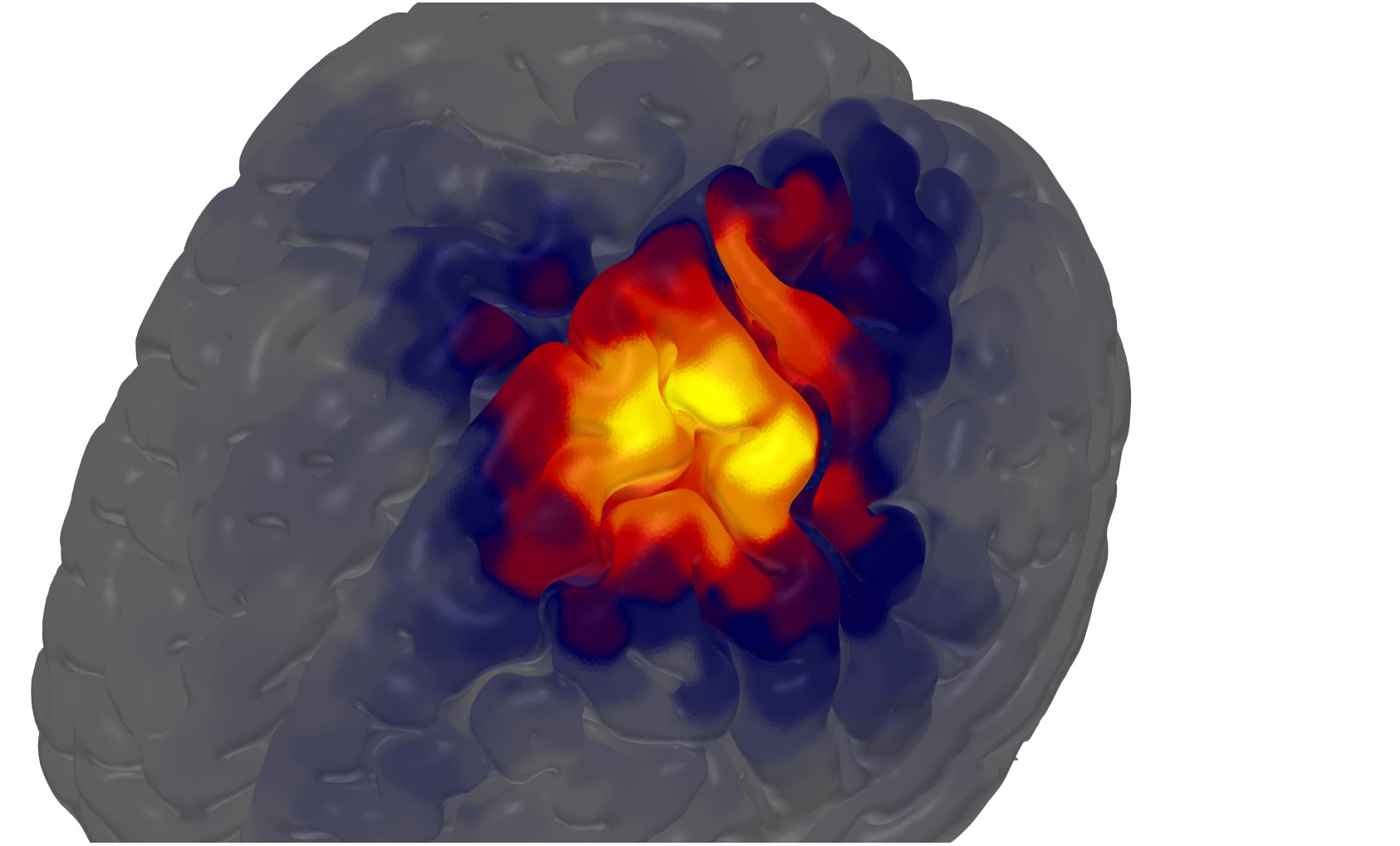}
        \end{minipage}%
        \begin{minipage}[b]{\imageWidth}
            \begin{center}
                Whitney difference
            \end{center}
            \includegraphics[trim={6cm 2cm 6.5cm 2cm},clip,width=1\linewidth]{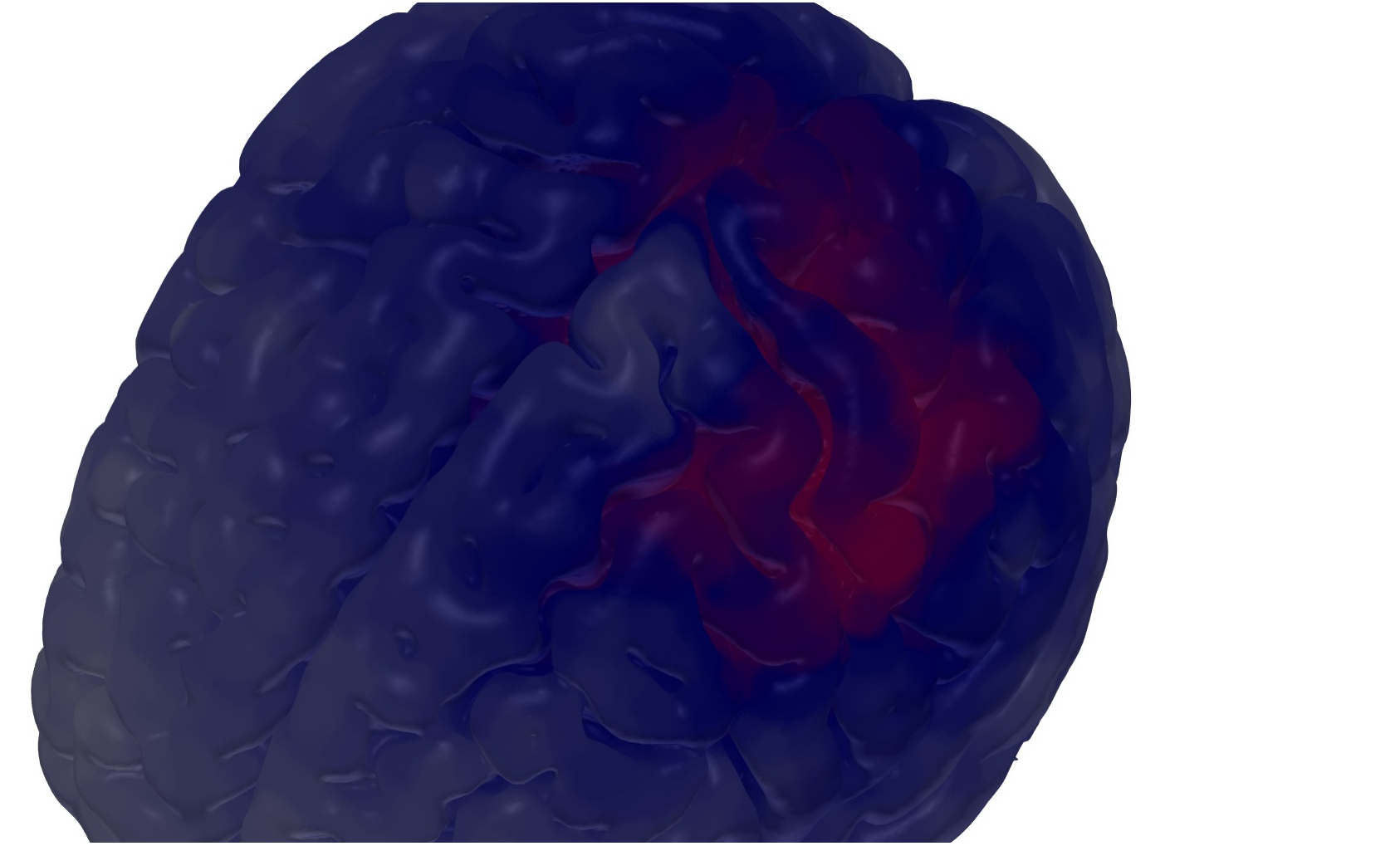}
        \end{minipage}%
        \begin{minipage}[b]{\imageWidth}
            \begin{center}
                $\Hdiv$ difference
            \end{center}
            \includegraphics[trim={6cm 2cm 6.5cm 2cm},clip,width=1\linewidth]{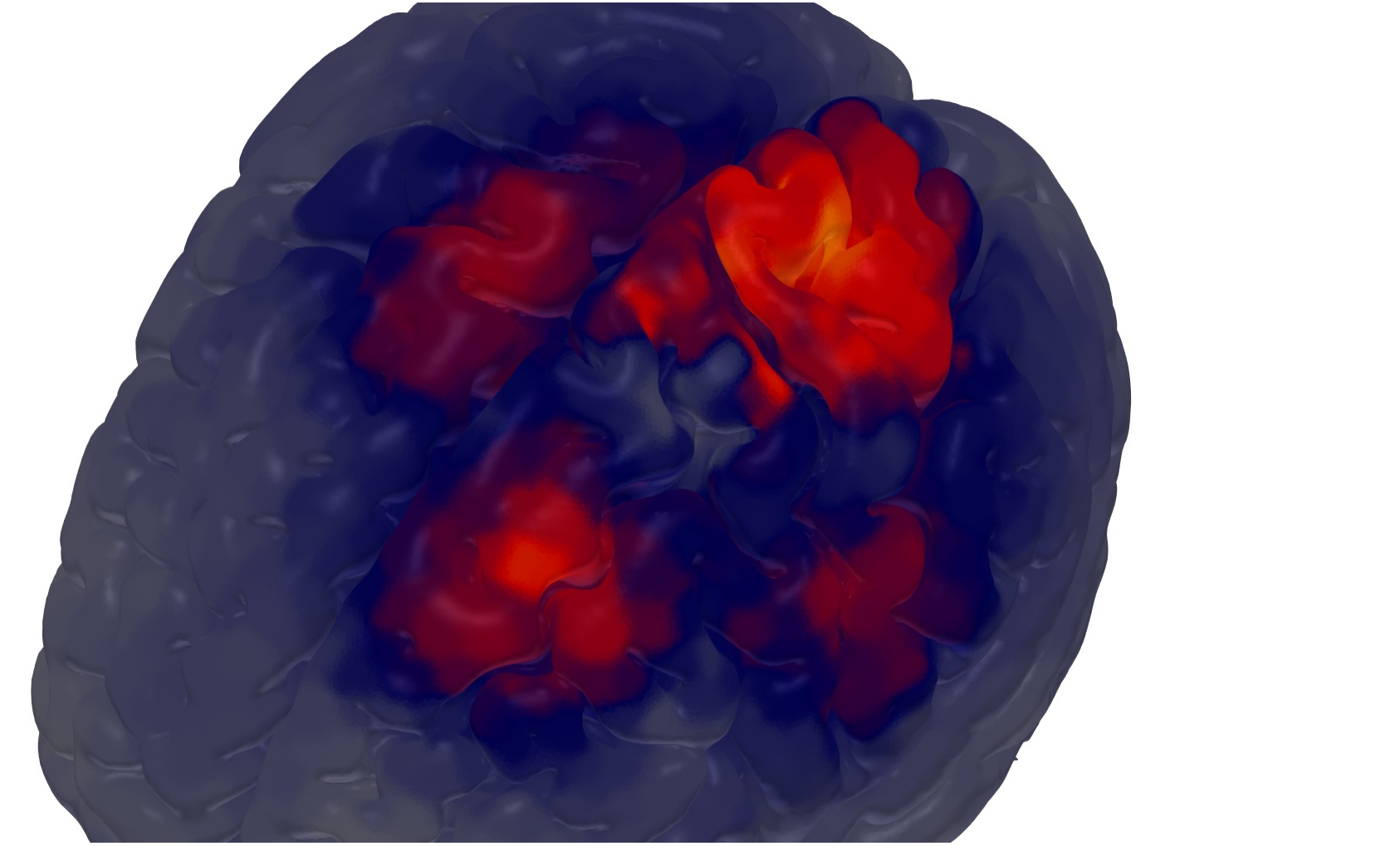}
        \end{minipage}

        \begin{minipage}[b]{\inverseMethodNameWidth}
            \rotatebox{90}{\hspace{0.05cm}SHAL1R}
        \end{minipage}%
        \begin{minipage}[b]{\imageWidth}
            \includegraphics[trim={6cm 2cm 6.5cm 2cm},clip,width=1\linewidth]{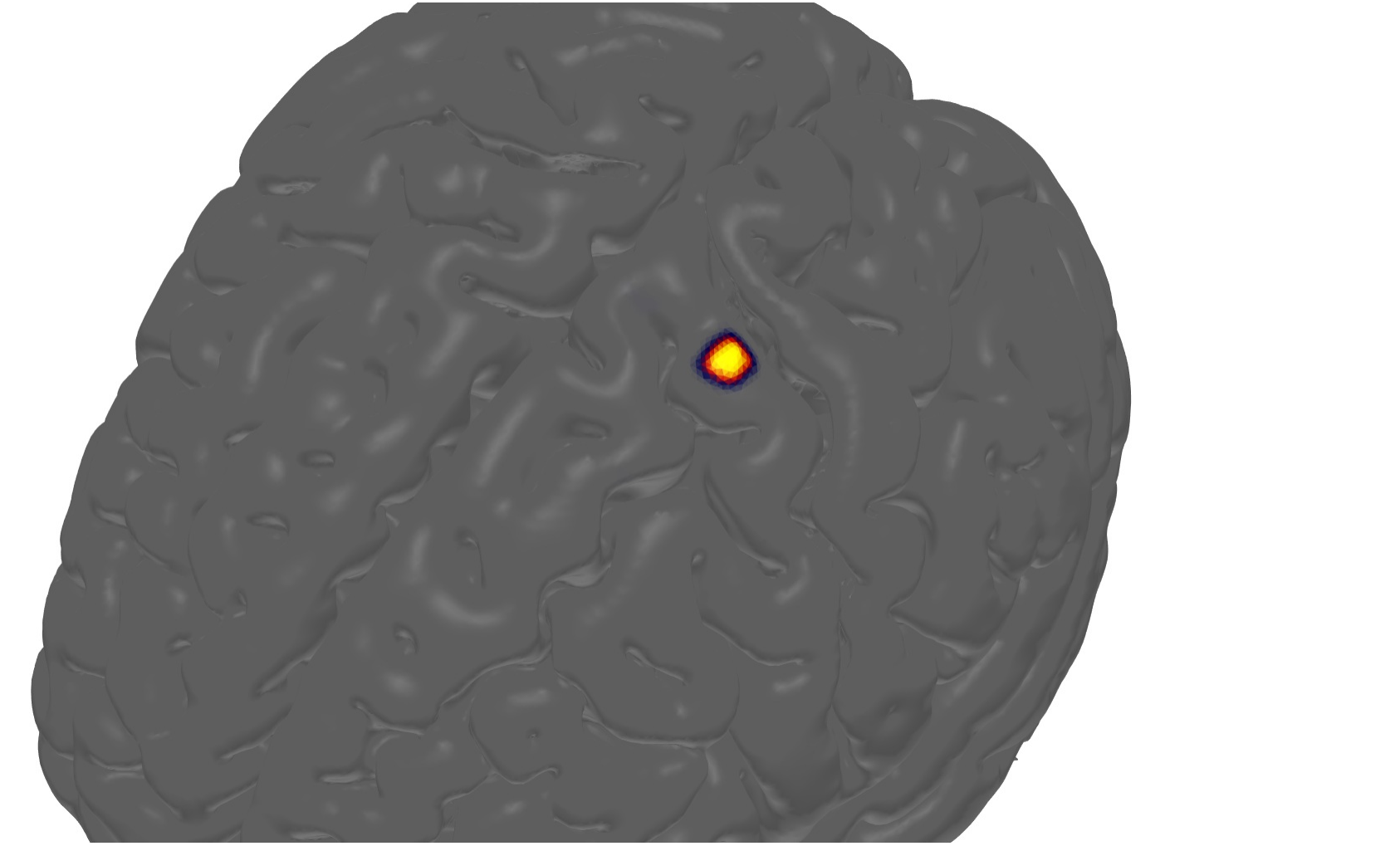}
        \end{minipage}%
        \begin{minipage}[b]{\imageWidth}
            \includegraphics[trim={6cm 2cm 6.5cm 2cm},clip,width=1\linewidth]{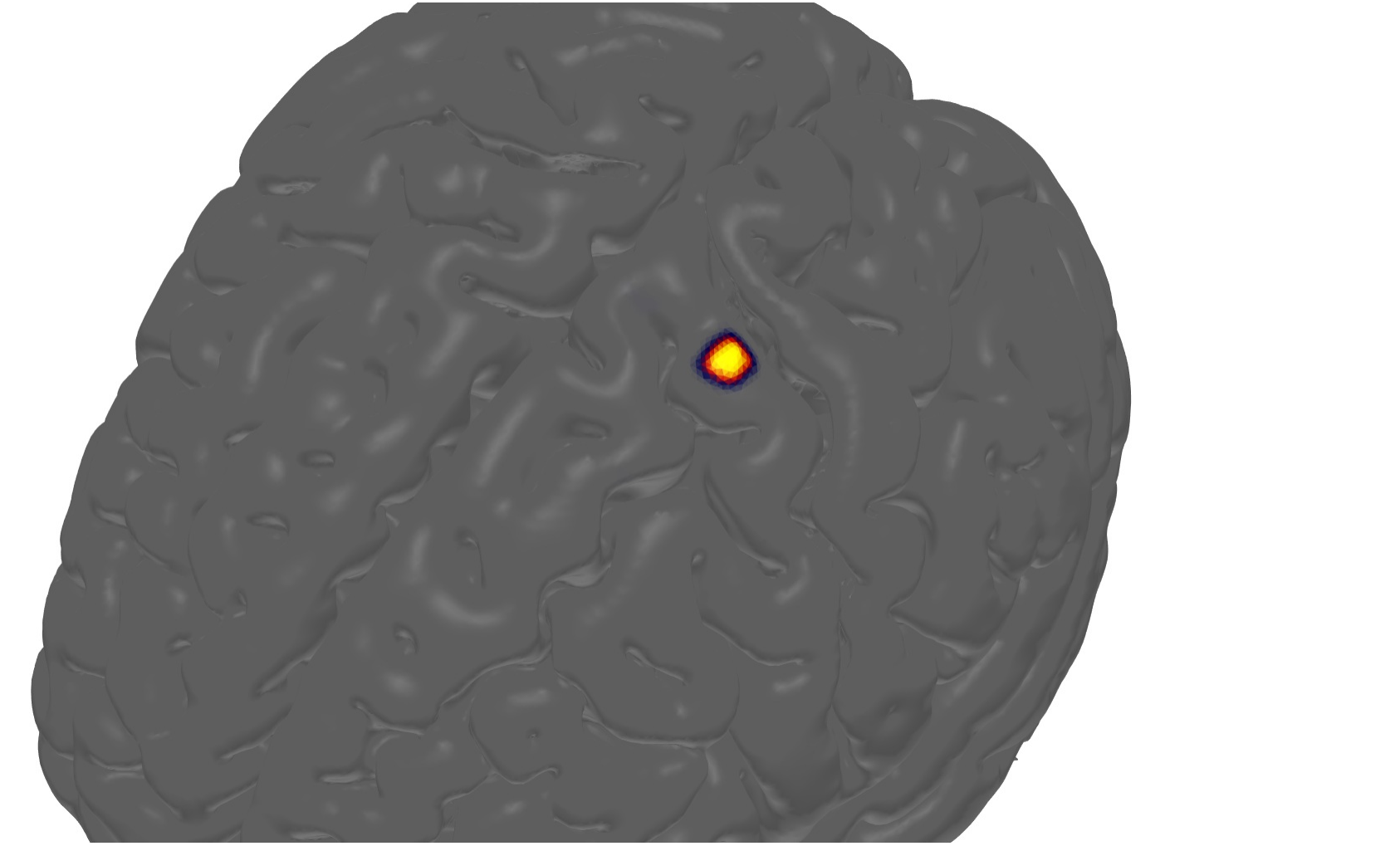}
        \end{minipage}%
        \begin{minipage}[b]{\imageWidth}
            \includegraphics[trim={6cm 2cm 6.5cm 2cm},clip,width=1\linewidth]{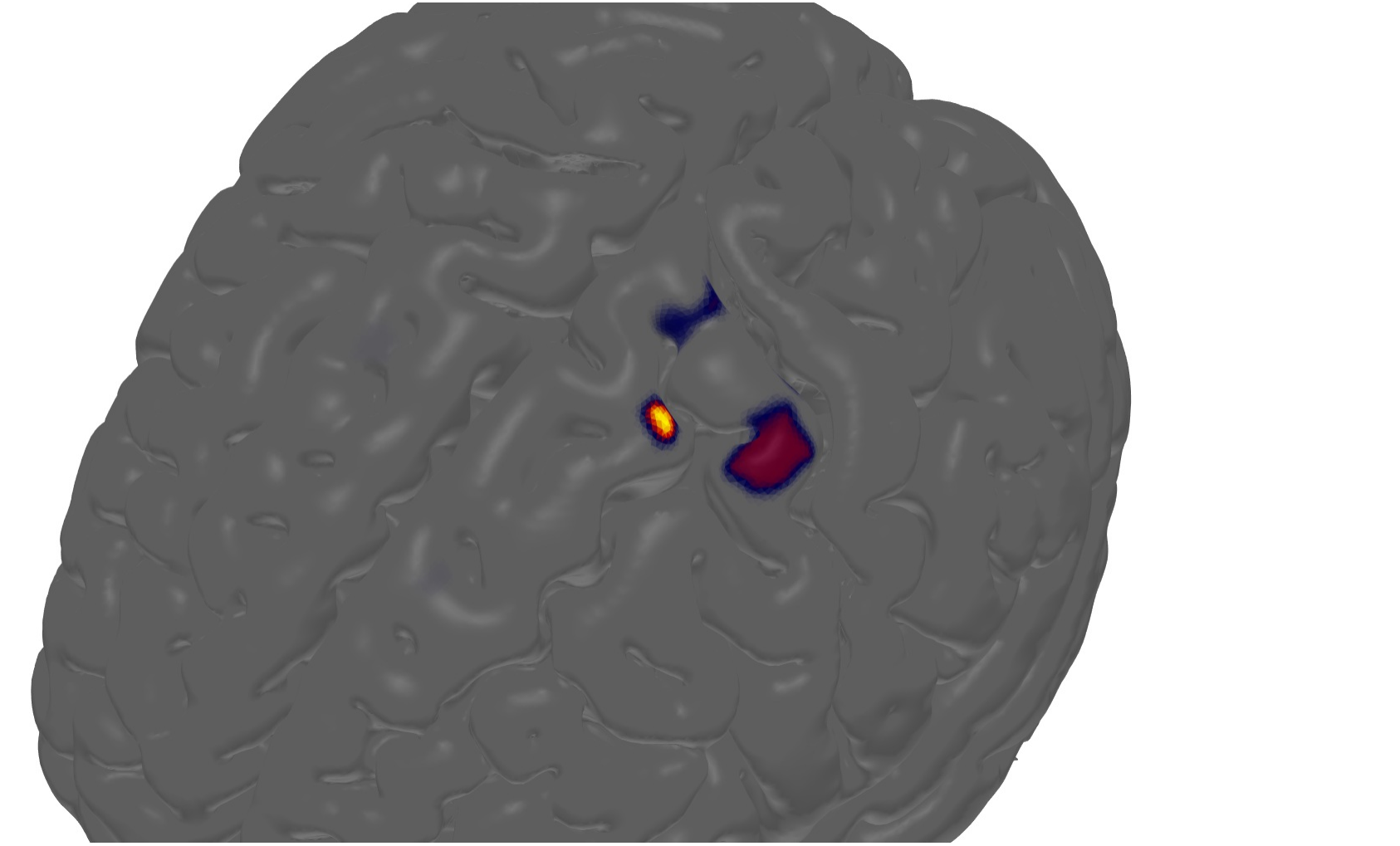}
        \end{minipage}%
        \begin{minipage}[b]{\imageWidth}
            \includegraphics[trim={6cm 2cm 6.5cm 2cm},clip,width=1\linewidth]{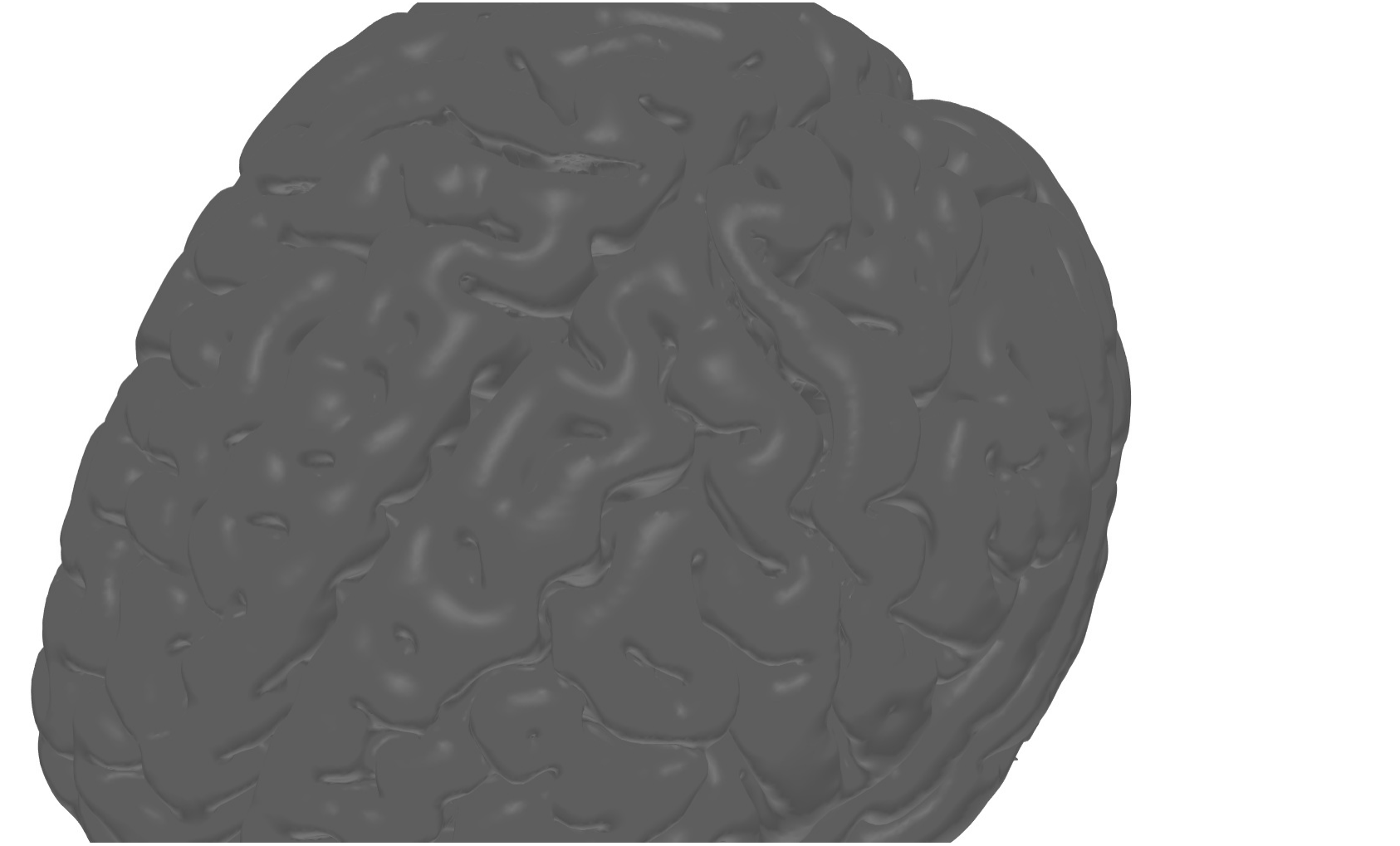}
        \end{minipage}%
        \begin{minipage}[b]{\imageWidth}
            \includegraphics[trim={6cm 2cm 6.5cm 2cm},clip,width=1\linewidth]{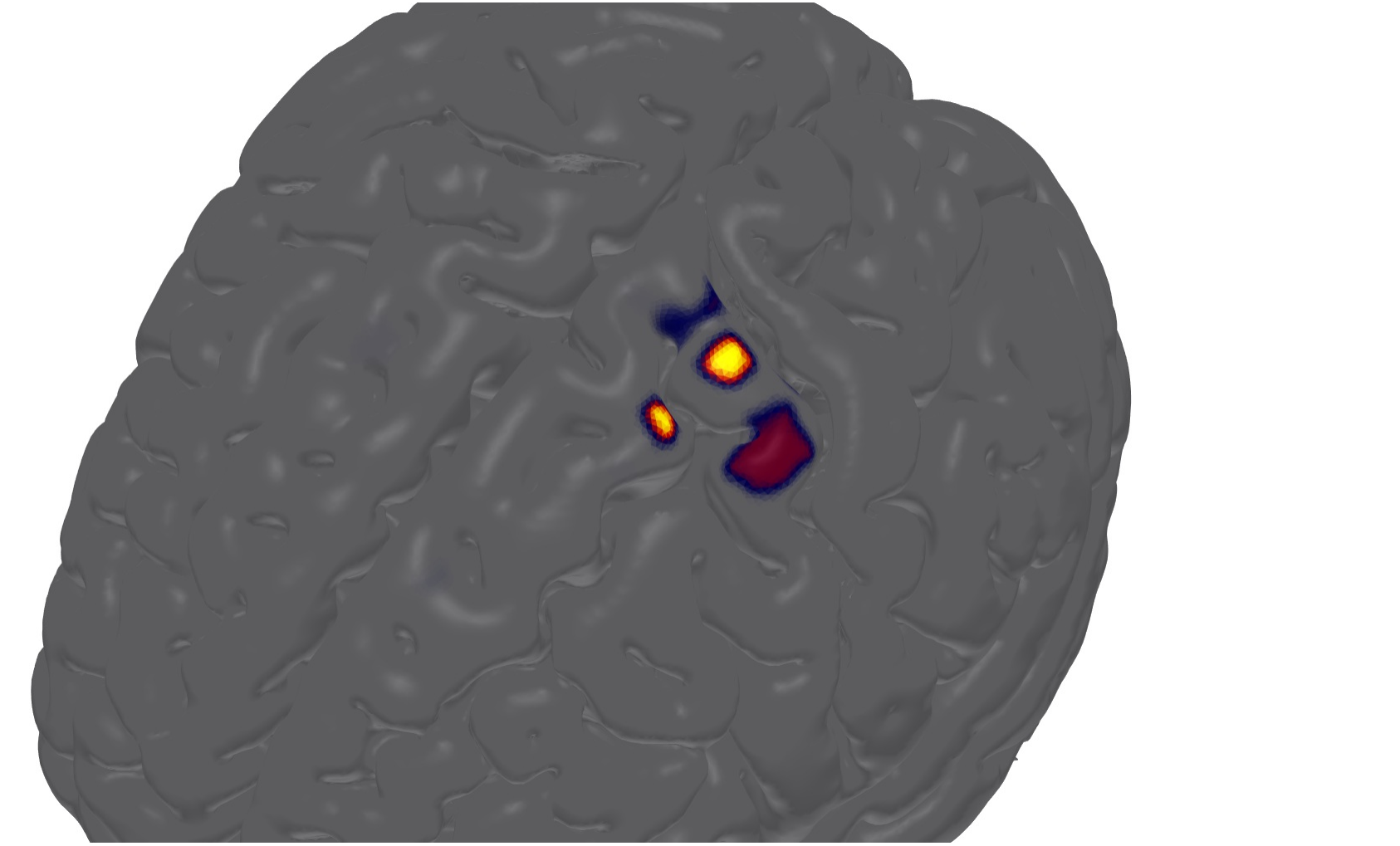}
        \end{minipage}

        \begin{minipage}[b]{\inverseMethodNameWidth}
            \rotatebox{90}{\hspace{0.2cm} SKF}
        \end{minipage}%
        \begin{minipage}[b]{\imageWidth}
            \includegraphics[trim={6cm 2cm 6.5cm 2cm},clip,width=1\linewidth]{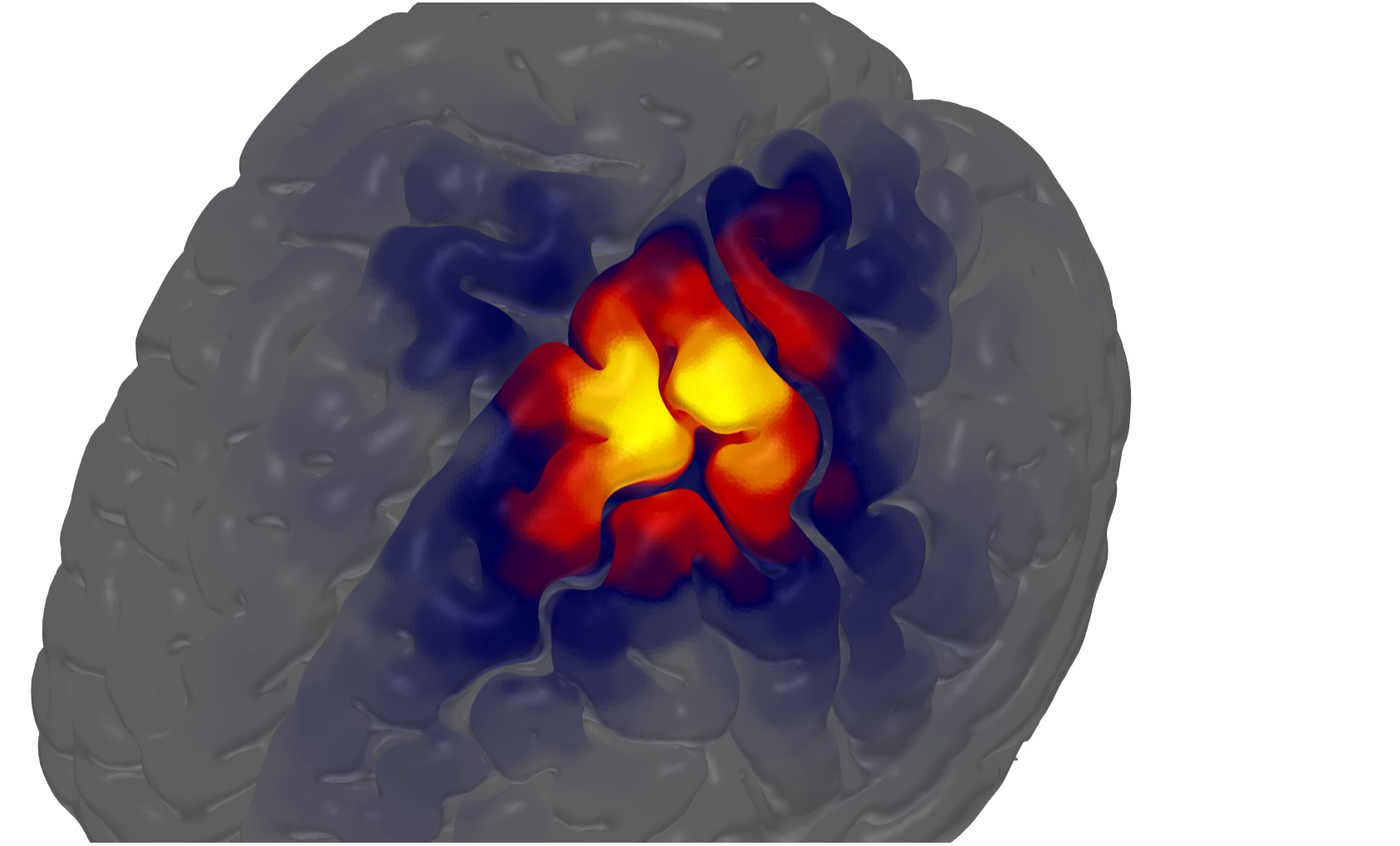}
        \end{minipage}%
        \begin{minipage}[b]{\imageWidth}
            \includegraphics[trim={6cm 2cm 6.5cm 2cm},clip,width=1\linewidth]{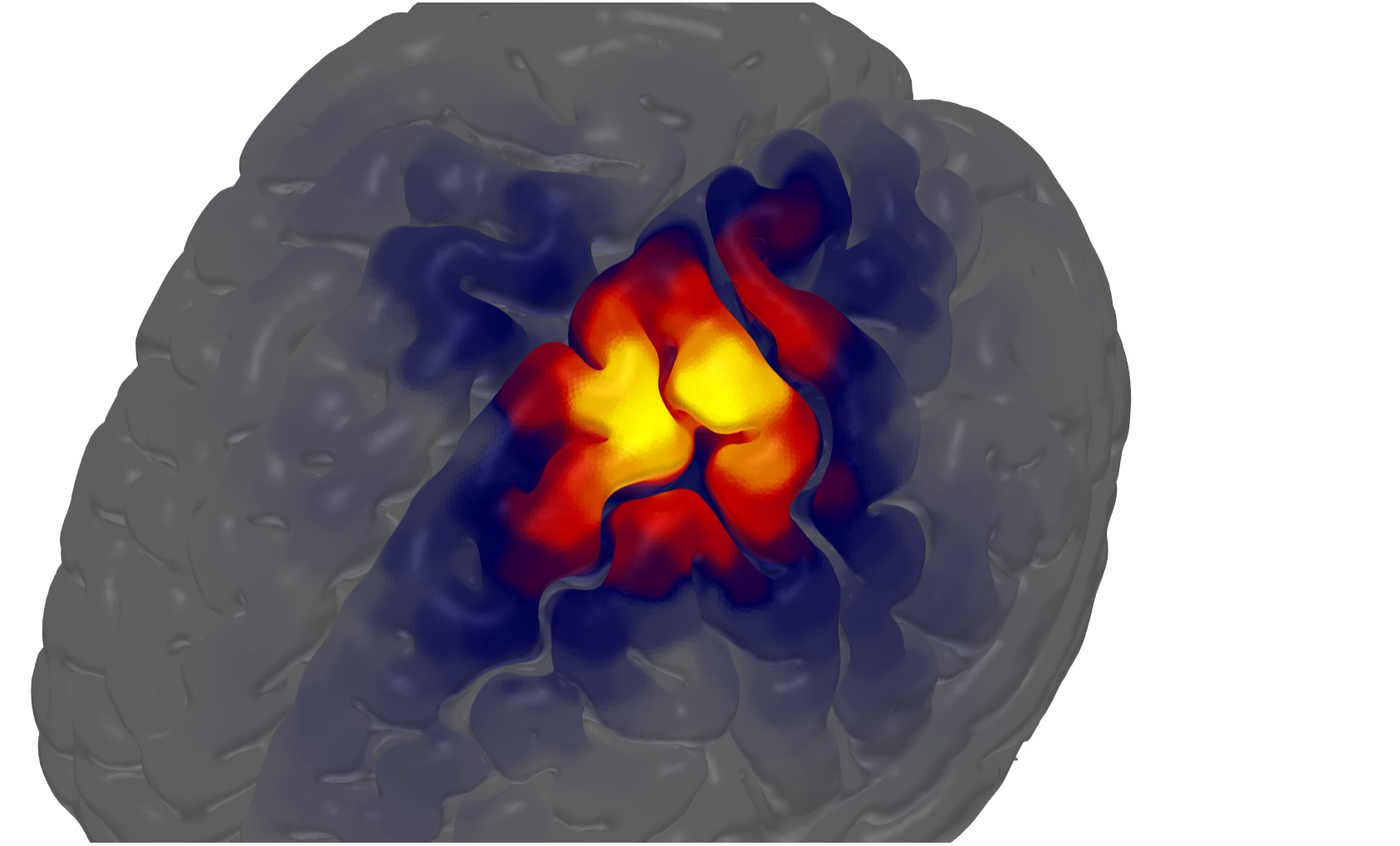}
        \end{minipage}%
        \begin{minipage}[b]{\imageWidth}
            \includegraphics[trim={6cm 2cm 6.5cm 2cm},clip,width=1\linewidth]{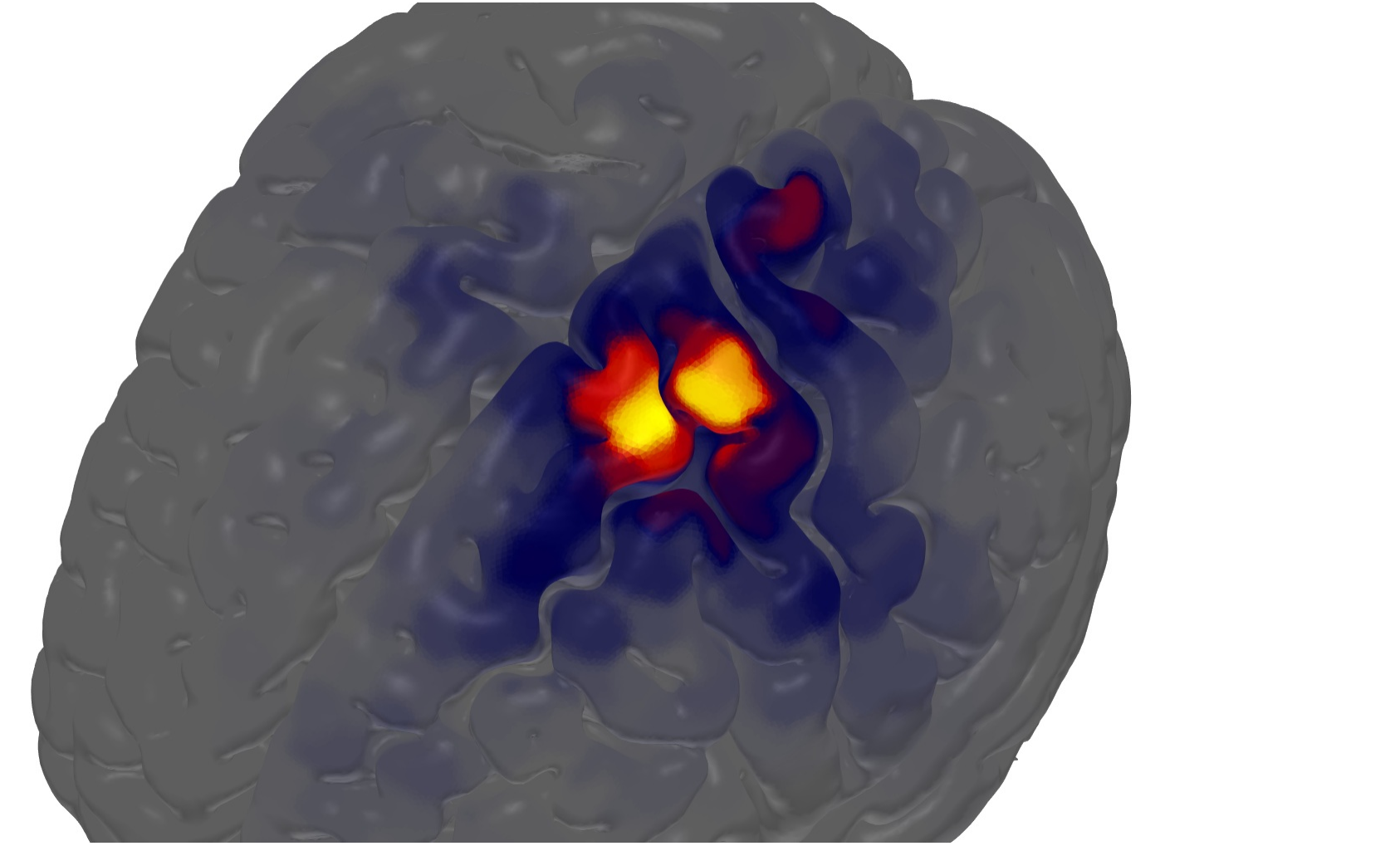}
        \end{minipage}%
        \begin{minipage}[b]{\imageWidth}
            \includegraphics[trim={6cm 2cm 6.5cm 2cm},clip,width=1\linewidth]{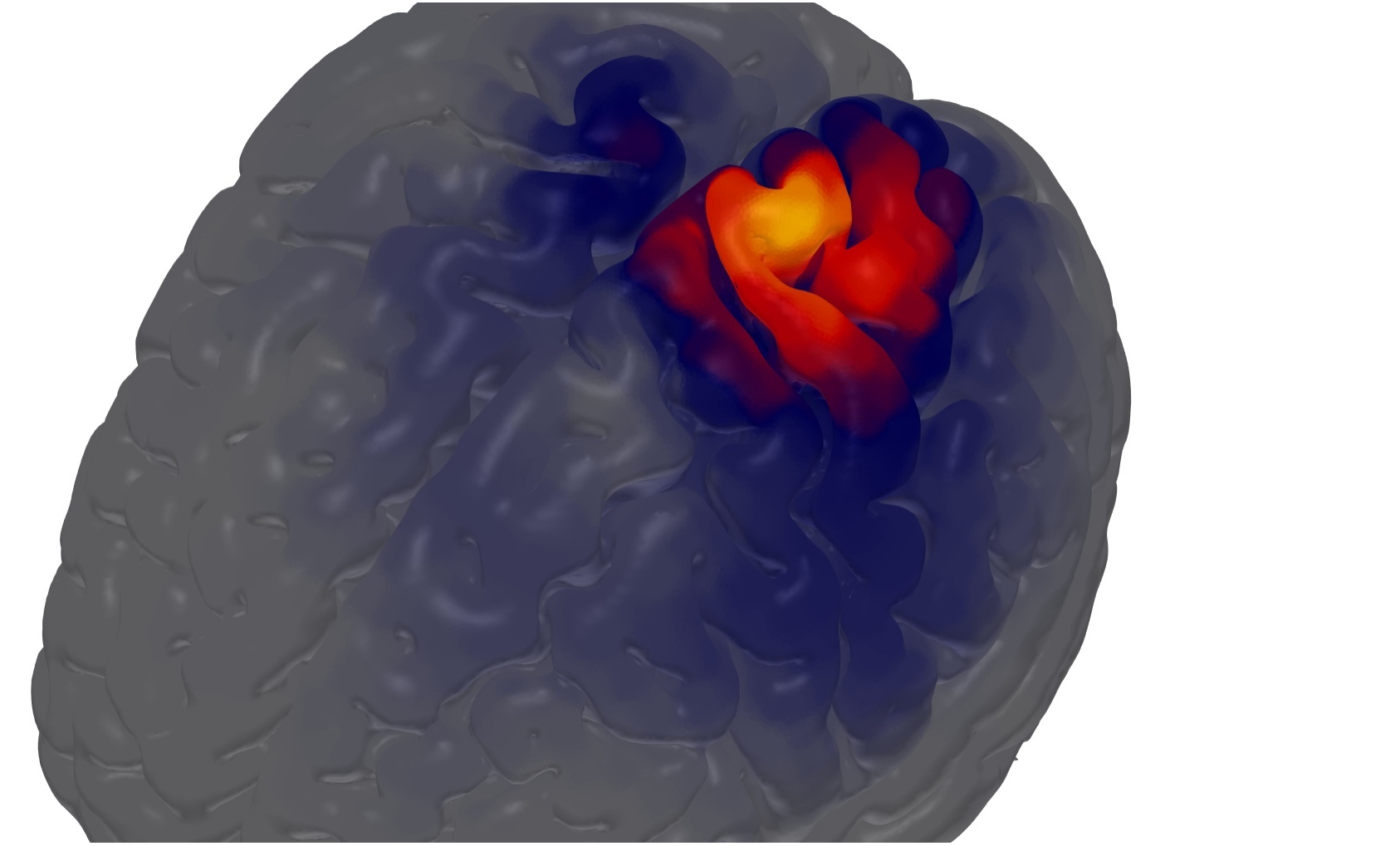}
        \end{minipage}%
        \begin{minipage}[b]{\imageWidth}
            \includegraphics[trim={6cm 2cm 6.5cm 2cm},clip,width=1\linewidth]{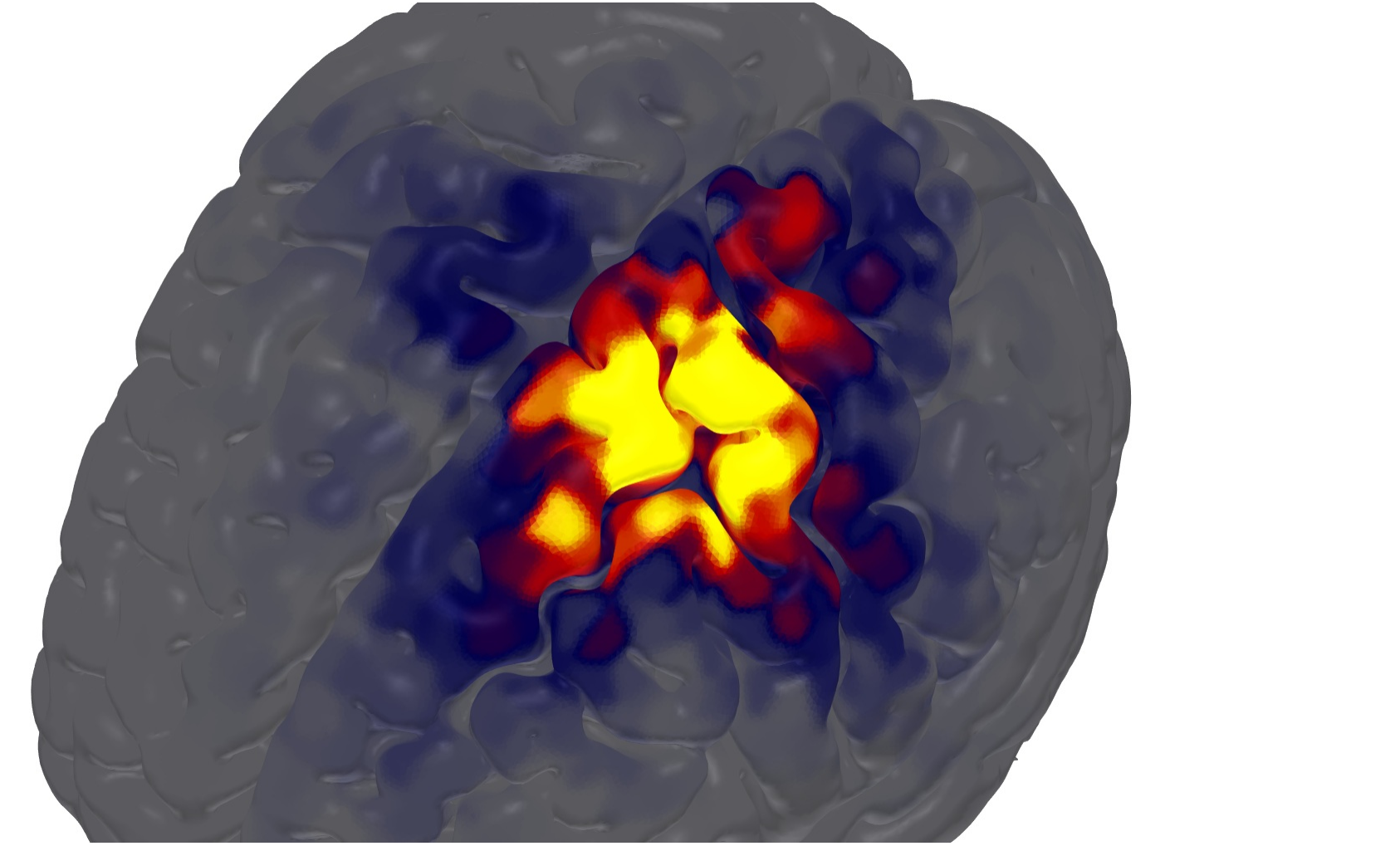}
        \end{minipage}

        \begin{minipage}[b]{\inverseMethodNameWidth}
            \rotatebox{90}{\hspace{0.3cm} DS}
        \end{minipage}%
        \begin{minipage}[b]{\imageWidth}
            \includegraphics[trim={6cm 2cm 6.5cm 2cm},clip,width=1\linewidth]{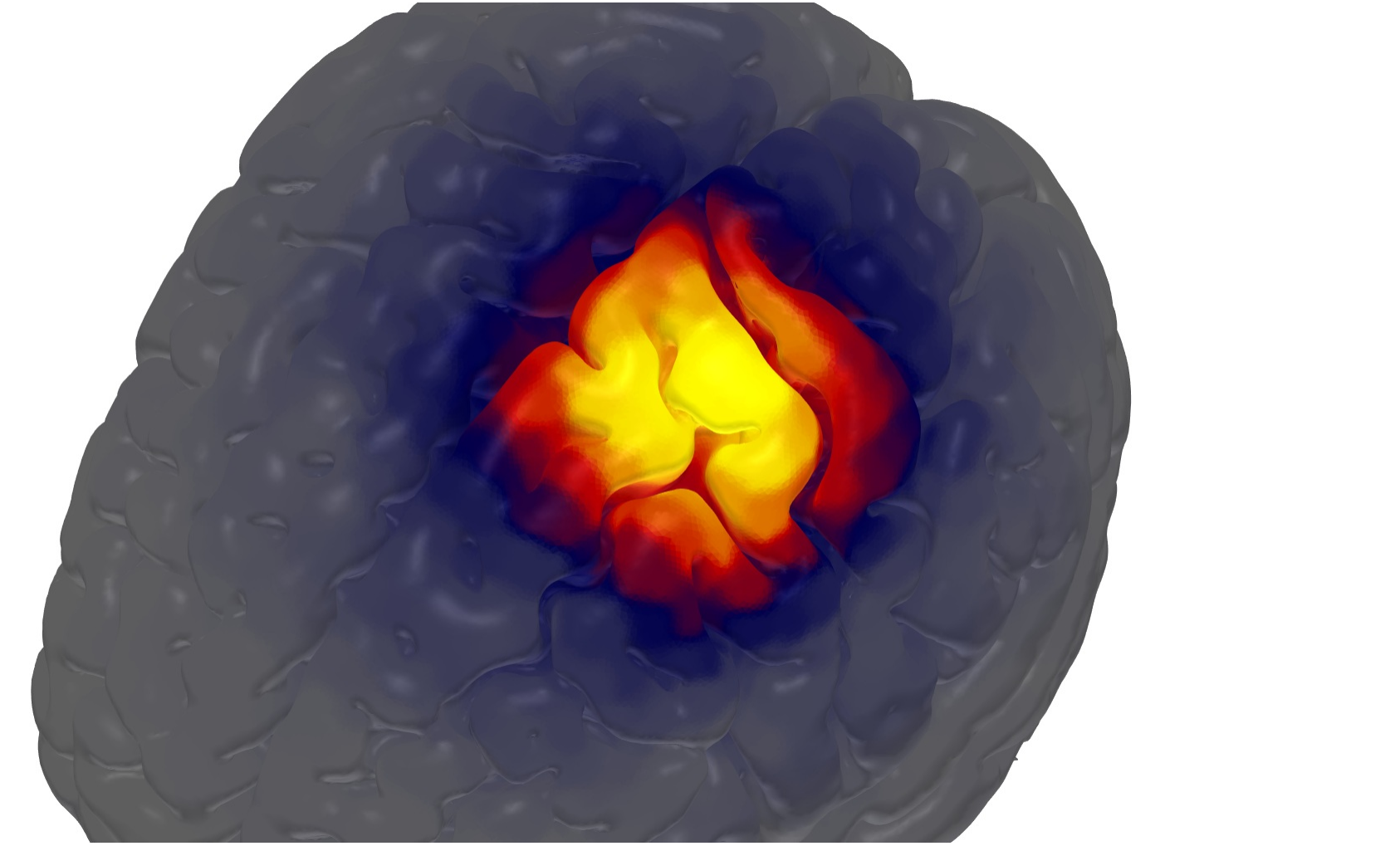}
        \end{minipage}%
        \begin{minipage}[b]{\imageWidth}
            \includegraphics[trim={6cm 2cm 6.5cm 2cm},clip,width=1\linewidth]{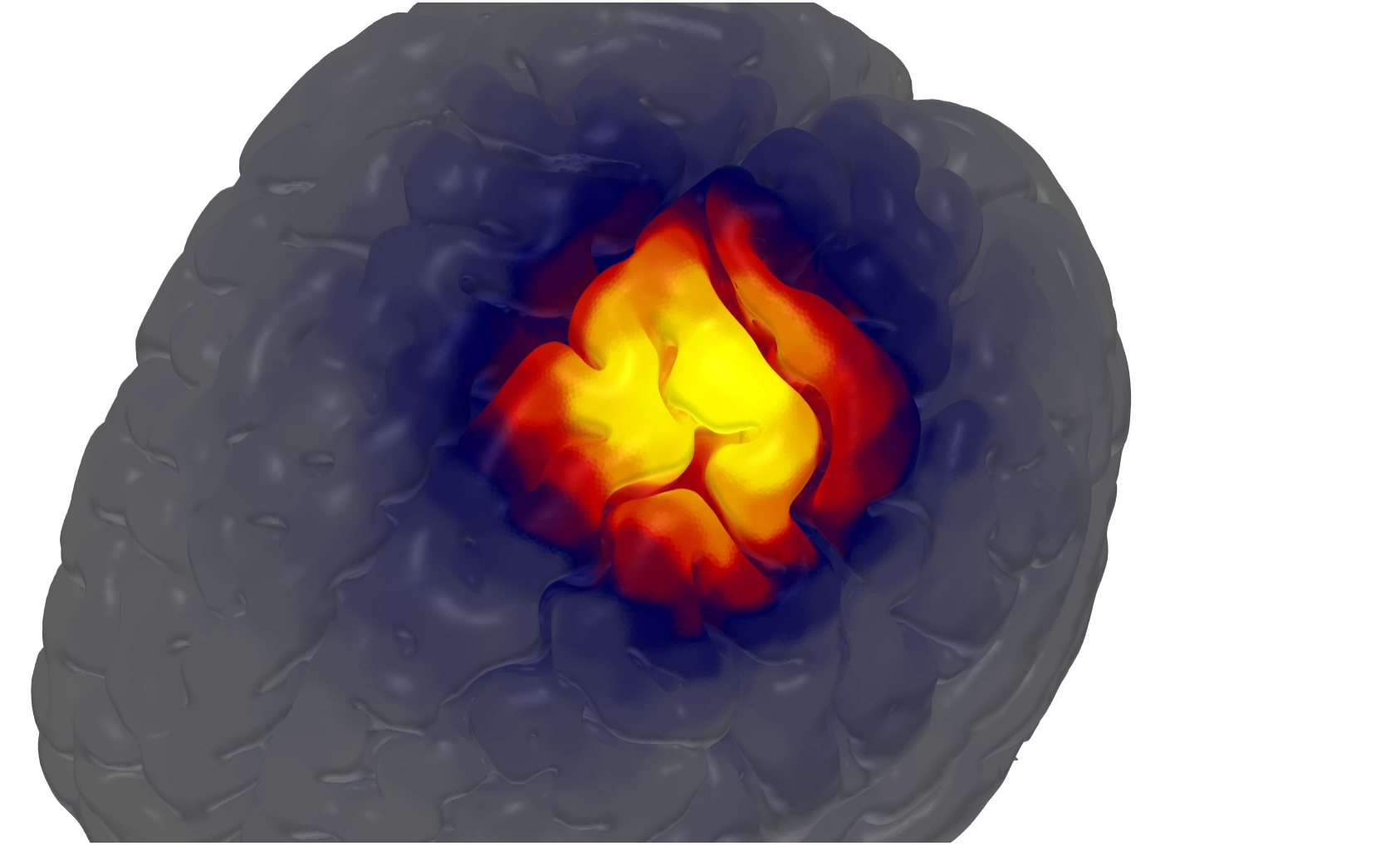}
        \end{minipage}%
        \begin{minipage}[b]{\imageWidth}
            \includegraphics[trim={6cm 2cm 6.5cm 2cm},clip,width=1\linewidth]{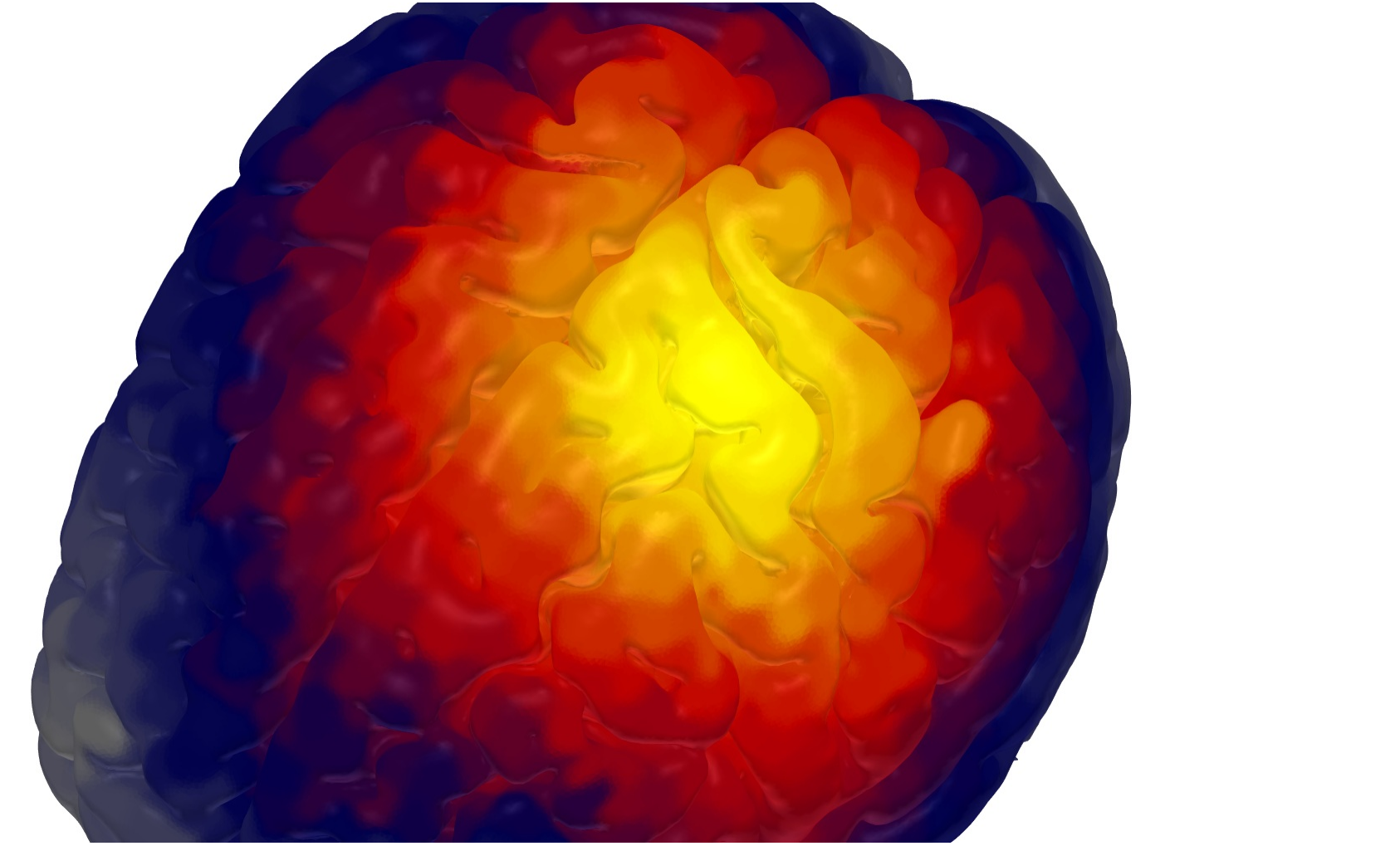}
        \end{minipage}%
        \begin{minipage}[b]{\imageWidth}
            \includegraphics[trim={6cm 2cm 6.5cm 2cm},clip,width=1\linewidth]{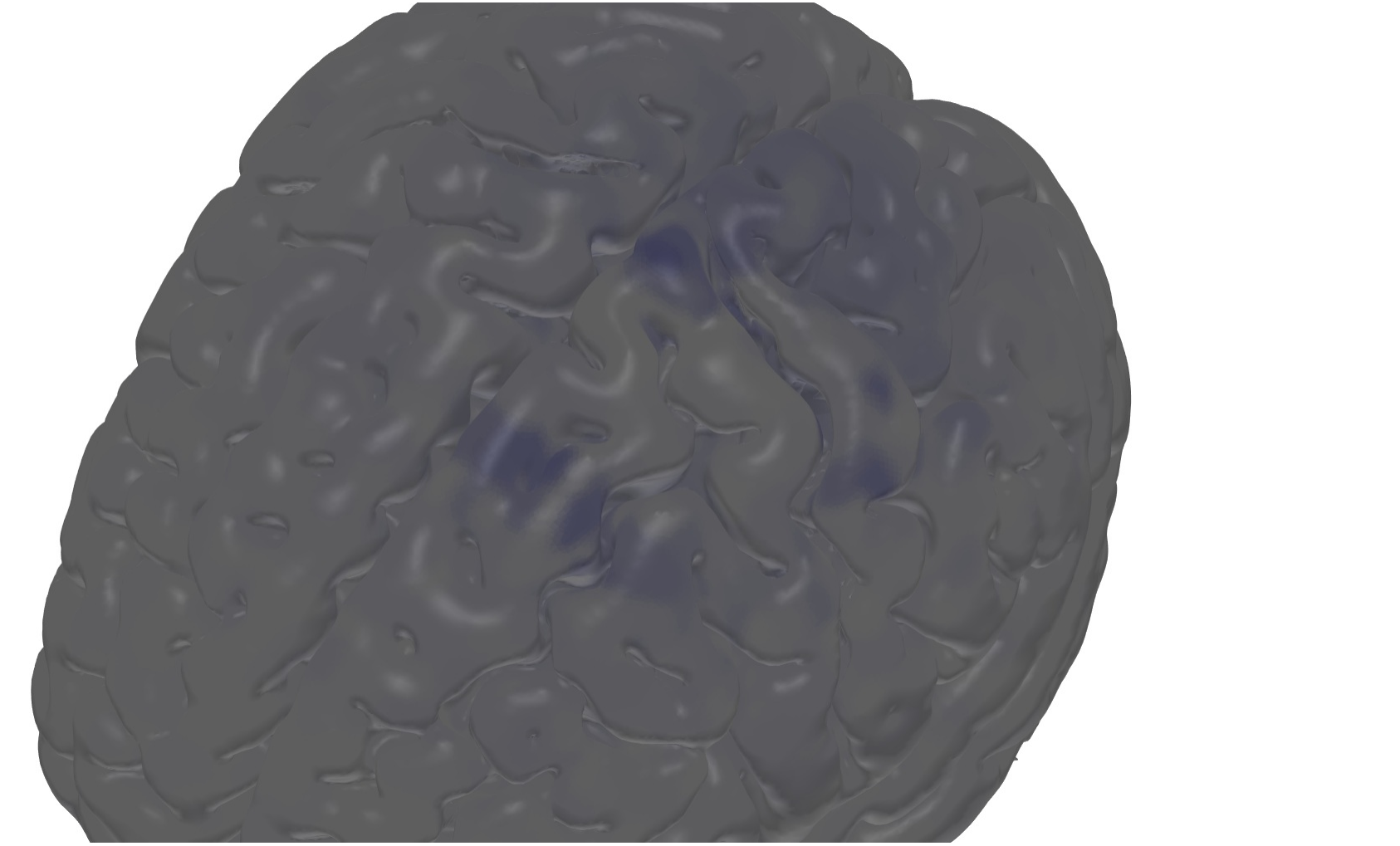}
        \end{minipage}%
        \begin{minipage}[b]{\imageWidth}
            \includegraphics[trim={6cm 2cm 6.5cm 2cm},clip,width=1\linewidth]{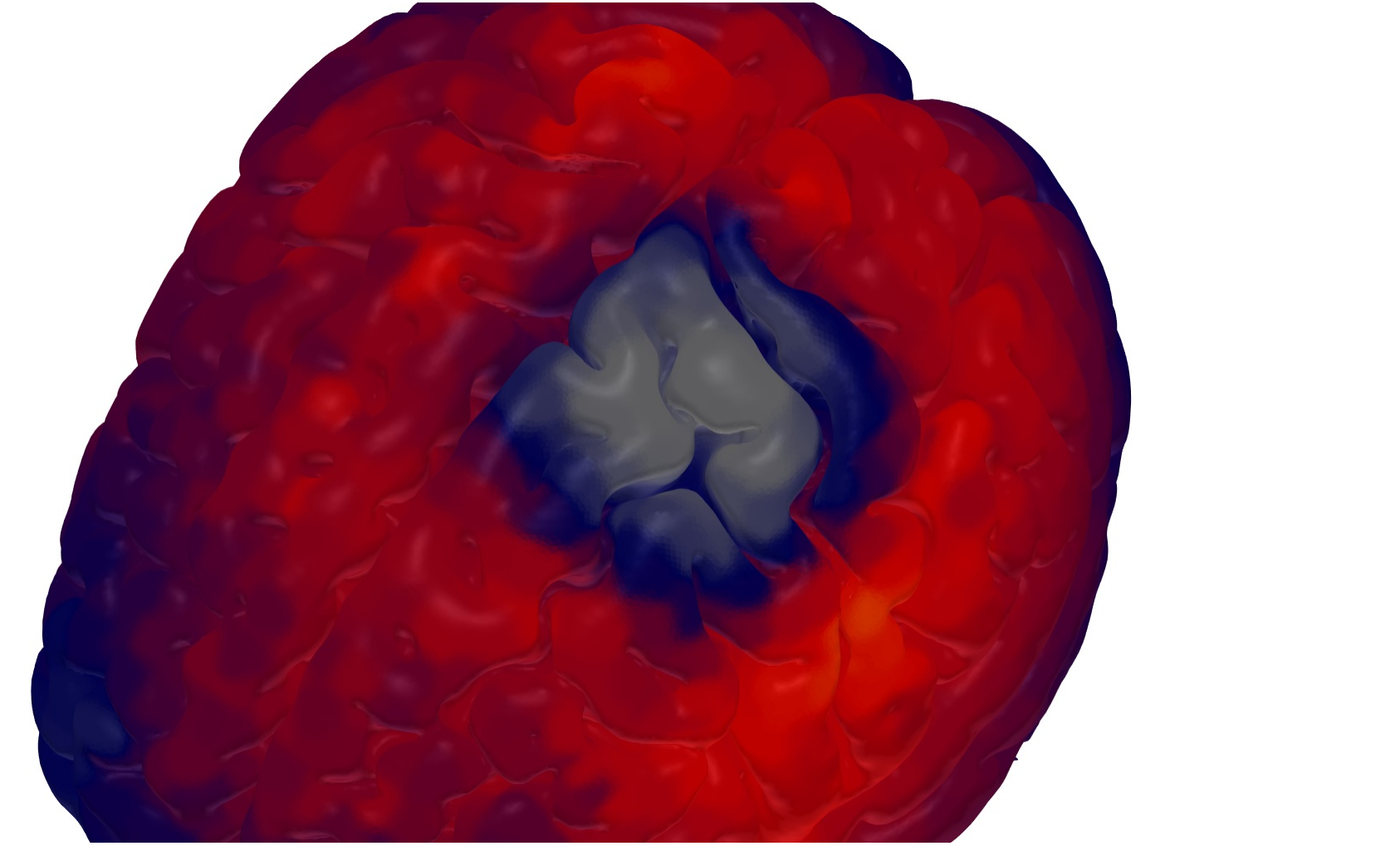}
        \end{minipage}
    \end{minipage}
    \hspace{0.01\linewidth}
    \begin{minipage}[b]{\colorbarWidth}
        \includegraphics[width=\linewidth]{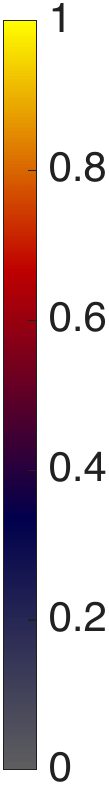}
    \end{minipage}
    \hfill
    \caption{\small Close up of reconstructions of a synthetic cortical dipole with sLORETA, SHAL1R, SKF and DS using DUNEuro's Whitney basis function implementation
    (column 1), DUNEuro's Local subtraction (column 2) and Zeffiro Interface's $\Hdiv$ (column 3). Column 4 shows the difference between Whitney and Local subtraction reconstructions, while column 5 does the same for $\Hdiv$ and Local subtraction. SHAL1R reconstructions show the true source location and do not differ from each other in the Whitney case. A SHAL1R reconstruction of a $\Hdiv$ source spreads the source across 3 positions near the true source position, resulting in an observable difference between it and a Local subtraction.}
    \label{fig.erotuskuva}
\end{figure}

Figure~\ref{fig.erotuskuva} displays how a reconstruction of a cortical source behaves when a lead field has been produced by interpolating $\transferMatrix$ via Whitney's face-intersecting and edgewise basis functions, Local subtraction or Zeffiro's $\Hdiv$, and then reconstructed with either sLORETA, SHAL1R, SKF, or DS. We also observe differences between Whitney and Local subtraction and $\Hdiv$ and Local subtraction.

It seems that SHAL1R is the most robust inverse method of the \num 4 presented here: not only is the reconstruction extremely focused around the actual single source position $\sourcePosition$, but the chosen lead field formulation does not have an effect on the reconstruction, as the difference between the Whitney and Local subtraction reconstructions is equally zero. sLORETA, SKF, and DS produce similar distributed reconstructions around $\sourcePosition$. SKF does exhibit slightly greater focality near the actual source position than sLORETA and DS, but at the same time, it is less robust to changes in the chosen lead-field interpolation approach, as exemplified by much larger relative differences between Local subtraction and both Whitney and $\Hdiv$ reconstructions. Out of these four, the smallest differences near the actual source position are obtained with DS, although after a certain distance from the true source SHAL1R reconstructions coincide better. The results obtained with the $\Hdiv$ source model and the reconstruction difference against Local subtraction show that $\Hdiv$ yields wider reconstructions for sLORETA, SHAL1R, and DS, as can be seen from the reconstruction differences. SKF, however, exhibits more focal estimation than the ones obtained with Local subtraction and the Whitney basis.

\begin{figure}[t]
    \centering
    \def\figWidth{0.45\linewidth}
    \begin{minipage}{0.05\linewidth}
        \rotatebox{90}{Whitney}
    \end{minipage}\begin{minipage}{\figWidth}
    \centering
    sLORETA \vspace{0.1cm}

        \includegraphics[width=0.98\linewidth]{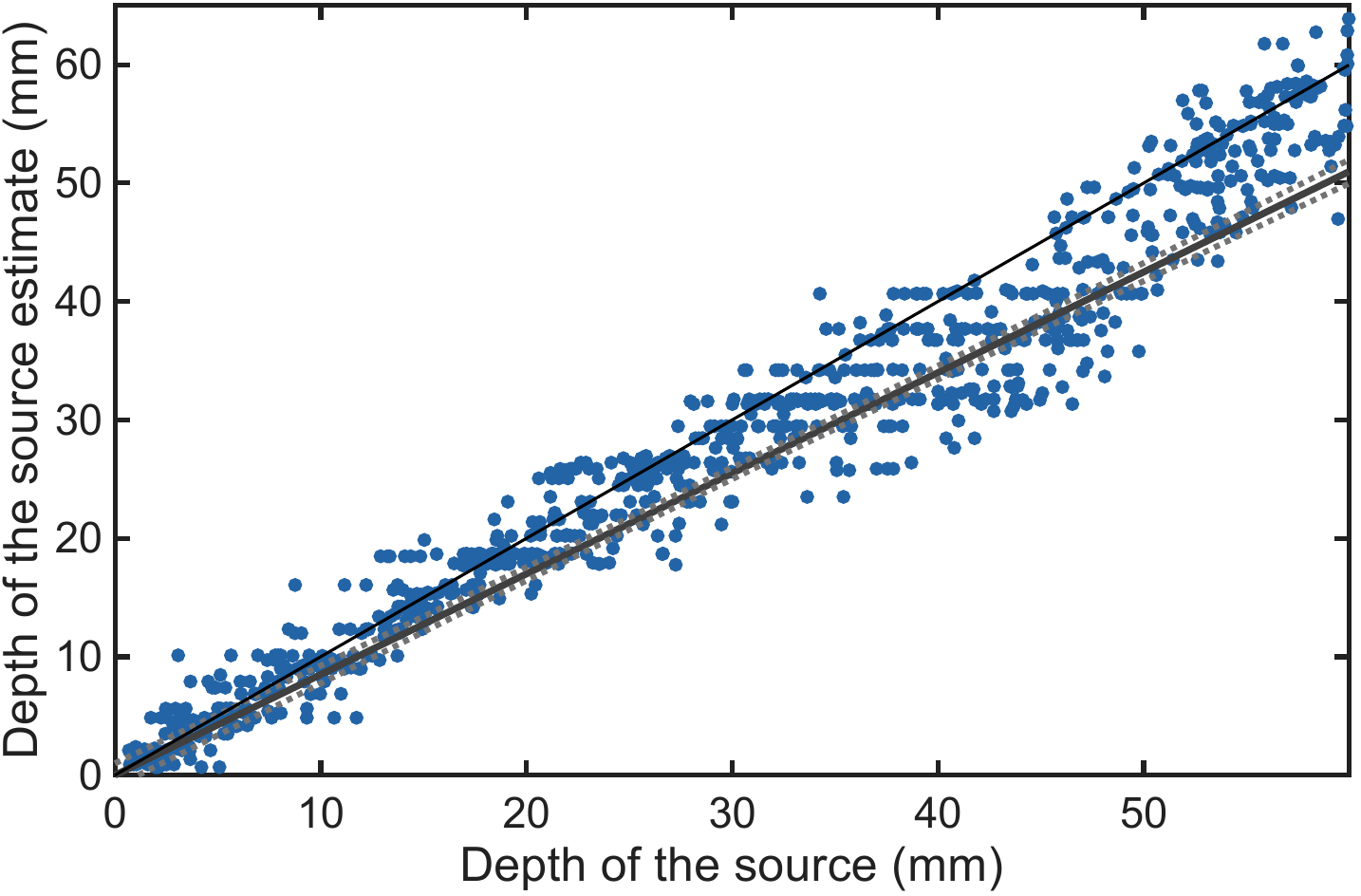}
    \end{minipage}%
    \begin{minipage}{\figWidth}
    \centering
    SHAL1R \vspace{0.1cm}

        \includegraphics[width=0.98\linewidth]{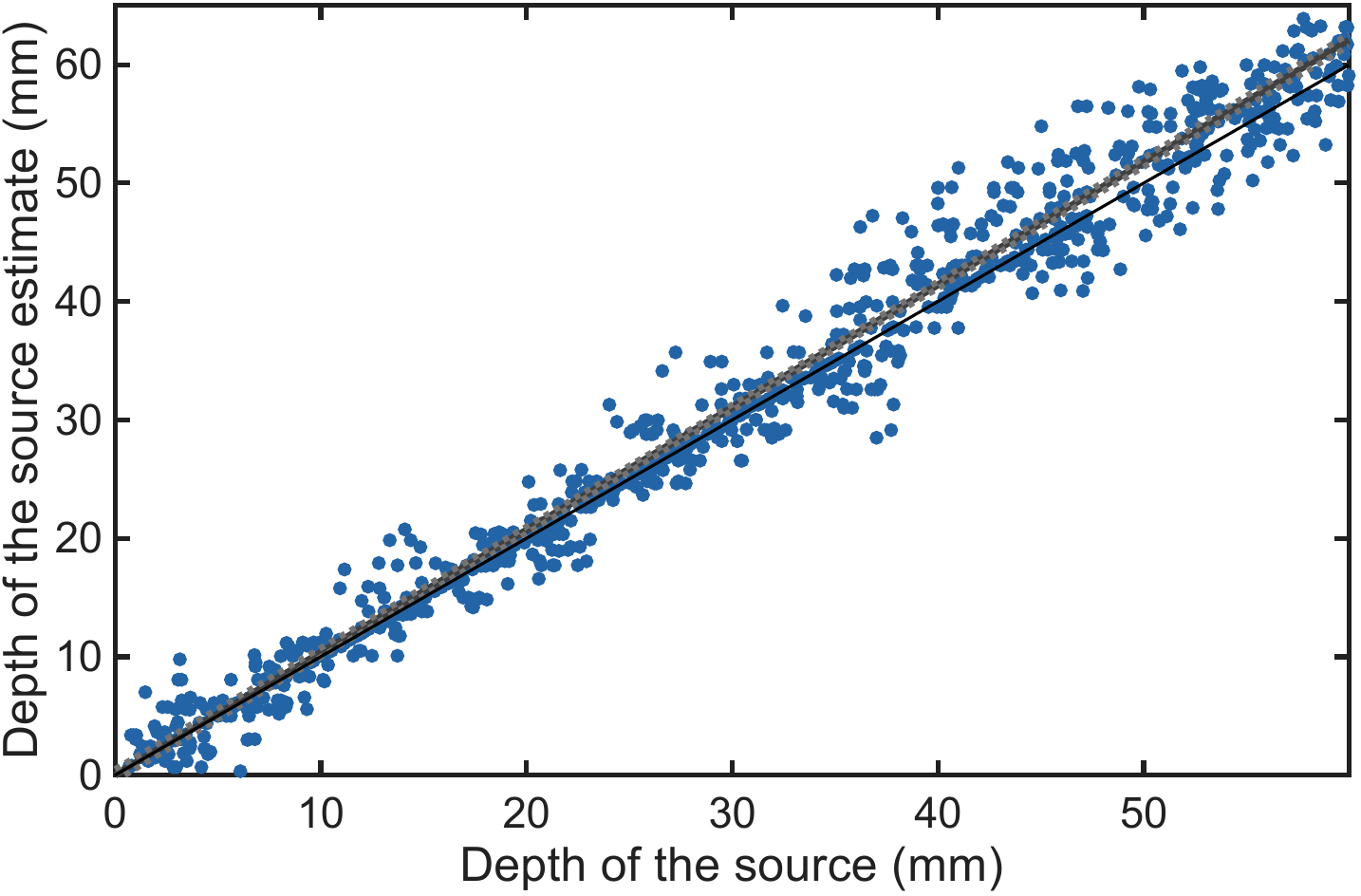}
    \end{minipage}

    \begin{minipage}{0.05\linewidth}
        \rotatebox{90}{$\Hdiv$}
    \end{minipage}%
    \begin{minipage}{\figWidth}
    \centering
        \includegraphics[width=0.98\linewidth]{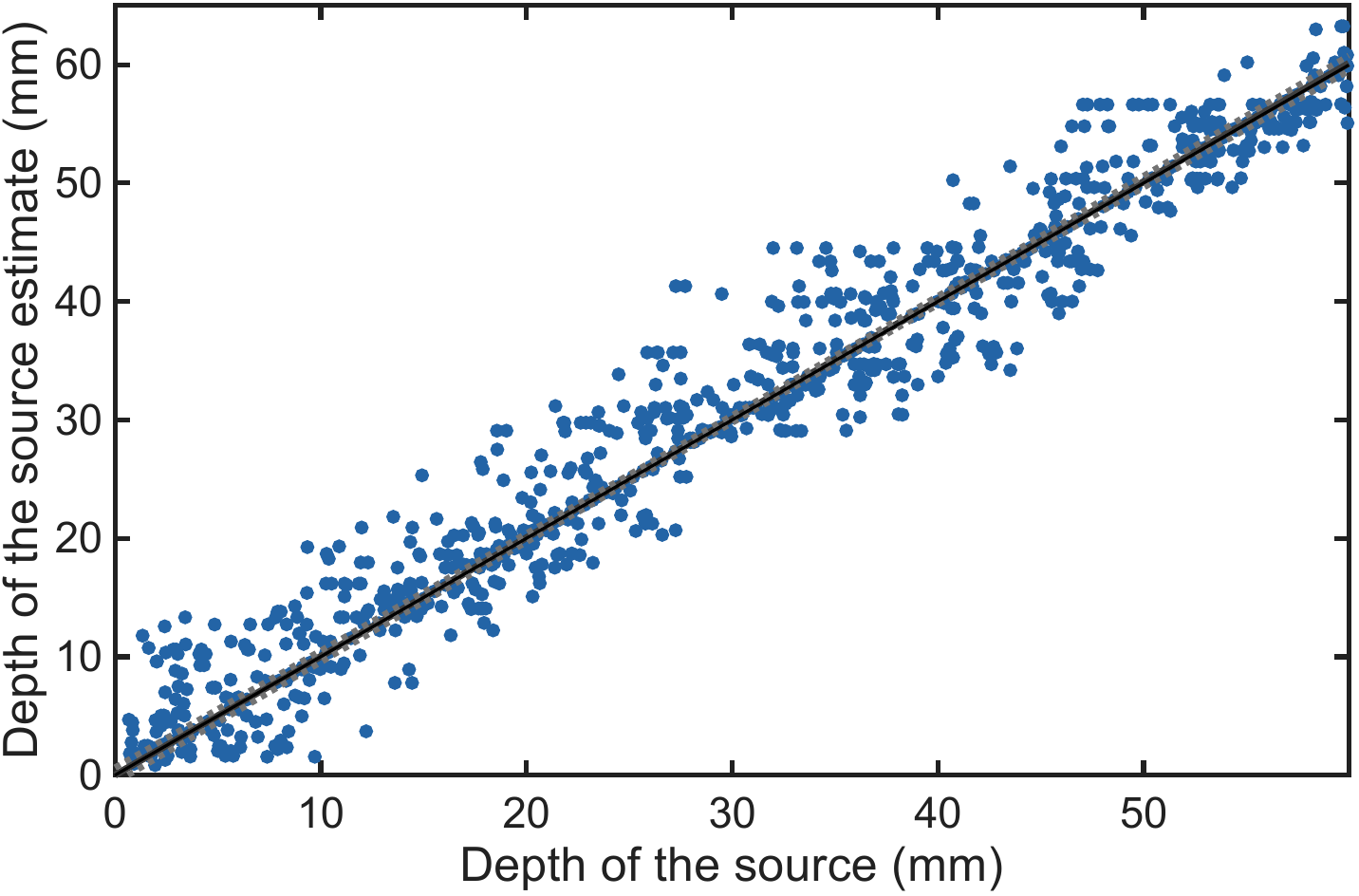}
    \end{minipage}\begin{minipage}{\figWidth}
    \centering
        \includegraphics[width=0.98\linewidth]{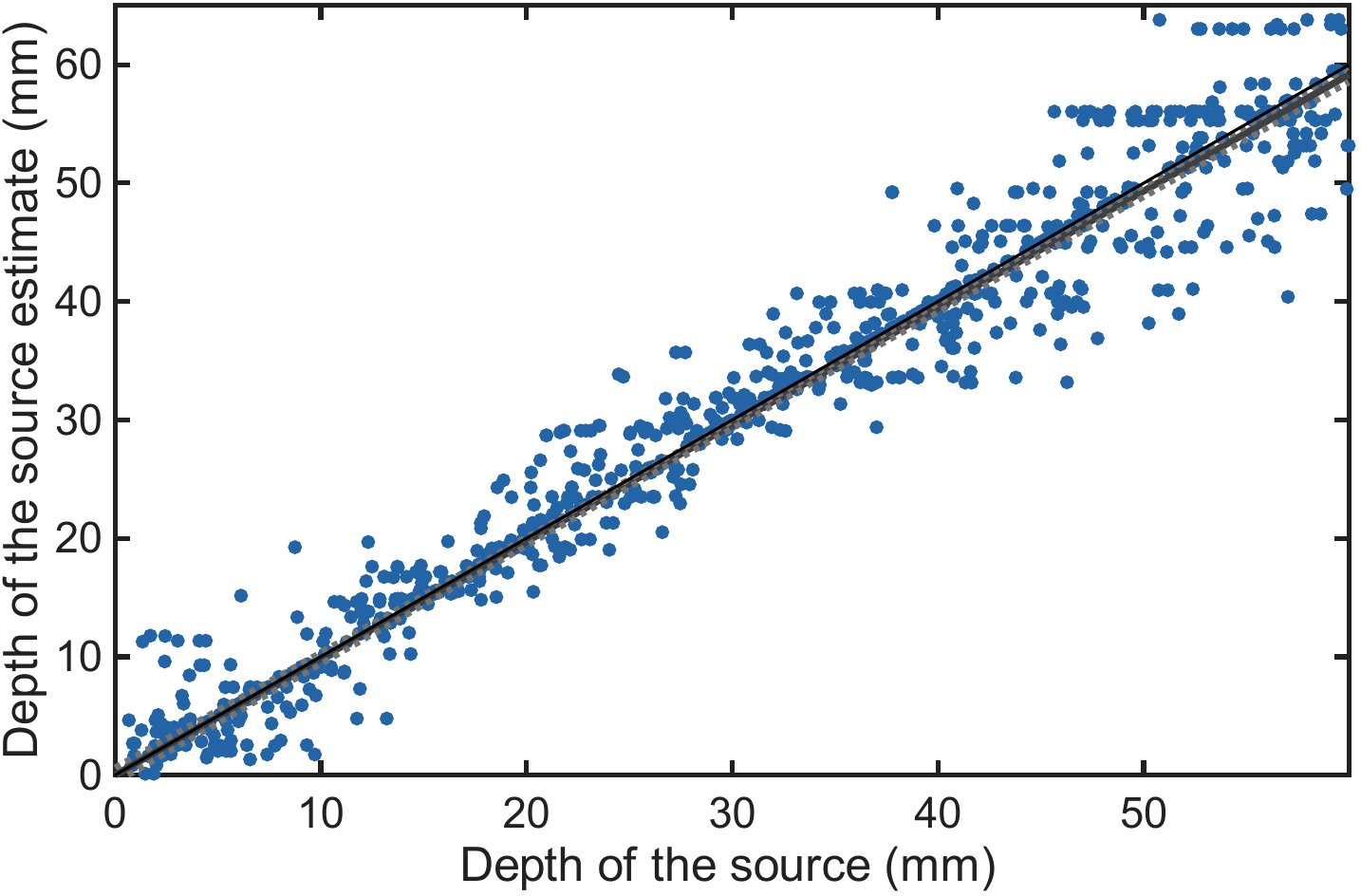}
    \end{minipage}

    \begin{minipage}{0.05\linewidth}
        \rotatebox{90}{Local subt.}
    \end{minipage}\begin{minipage}{\figWidth}
    \centering
        \includegraphics[width=0.98\linewidth]{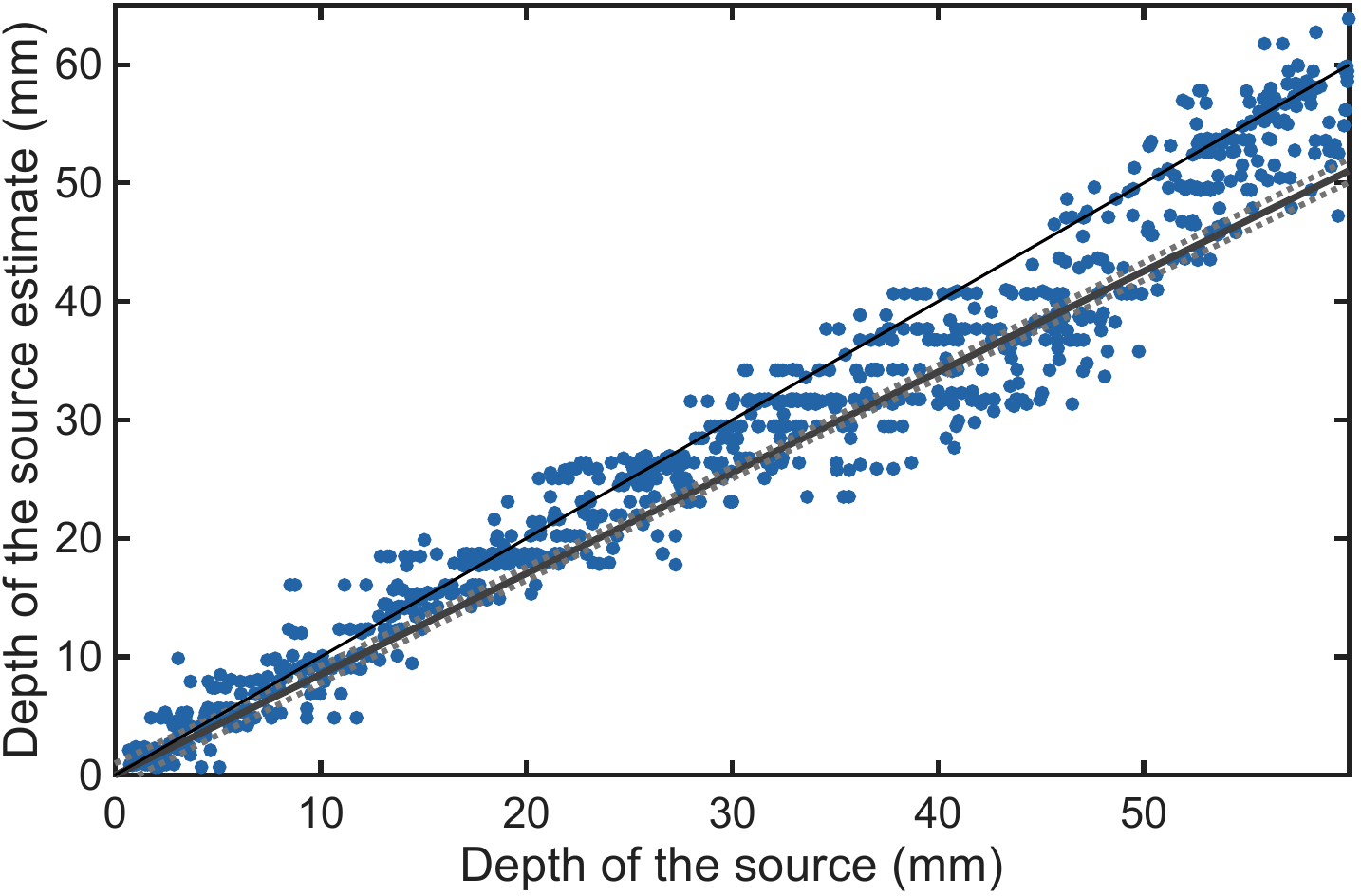}
    \end{minipage}%
    \begin{minipage}{\figWidth}
    \centering
        \includegraphics[width=0.98\linewidth]{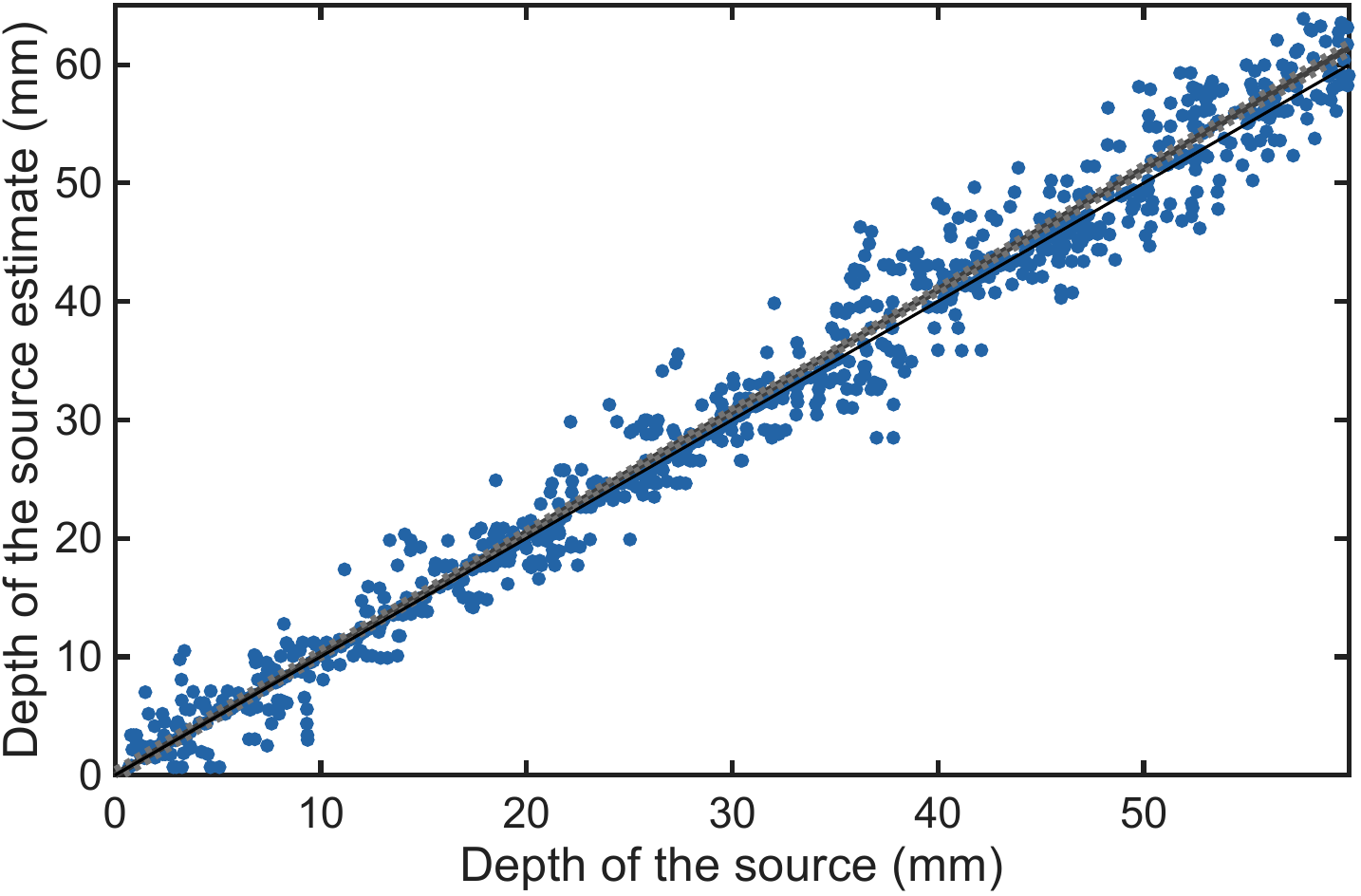}
    \end{minipage}

    \caption{\small Depth of true source plotted against the depth of estimation done by sLORETA and SHAL1R. The thin black line shows the optimal agreement between the true and estimated depth when the localization error is zero. The dark gray solid line displays the linear regression, and the dashed gray curves are the 95 \% confidence intervals.}
    \label{fig.depthbias}
\end{figure}

In Figure~\ref{fig.depthbias} we also observe how well sLORETA and SHAL1R reconstruct the depth or distance from the inner surface of the skull of a given source, when the interpolation schemes implied by Whitney basis functions, Local subtraction, and $\Hdiv$ are used.

It is obvious that SHAL1R agrees better with the optimal localization depth than sLORETA, as the slopes of the regression lines are 1.03 for Whitney and 1.02 for Local subtraction, whereas sLORETA yields 0.85 for both source models. In both cases, Local subtraction, as a forward interpolation scheme, produces slightly better localization results than Whitney, which again shows a better correspondence to the source--estimate regression line. The best regression lines are obtained with the $\Hdiv$ source model, since the sLORETA slope is 1.00 and the SHAL1R slope is 0.99.

\begin{figure}[t]
    \def\figWidth{0.45\linewidth}
    \centering
    \begin{minipage}{0.05\linewidth}
        \rotatebox{90}{Whitney}
    \end{minipage}\begin{minipage}{\figWidth}
    \centering
    sLORETA \vspace{0.1cm}

        \includegraphics[width=0.98\linewidth]{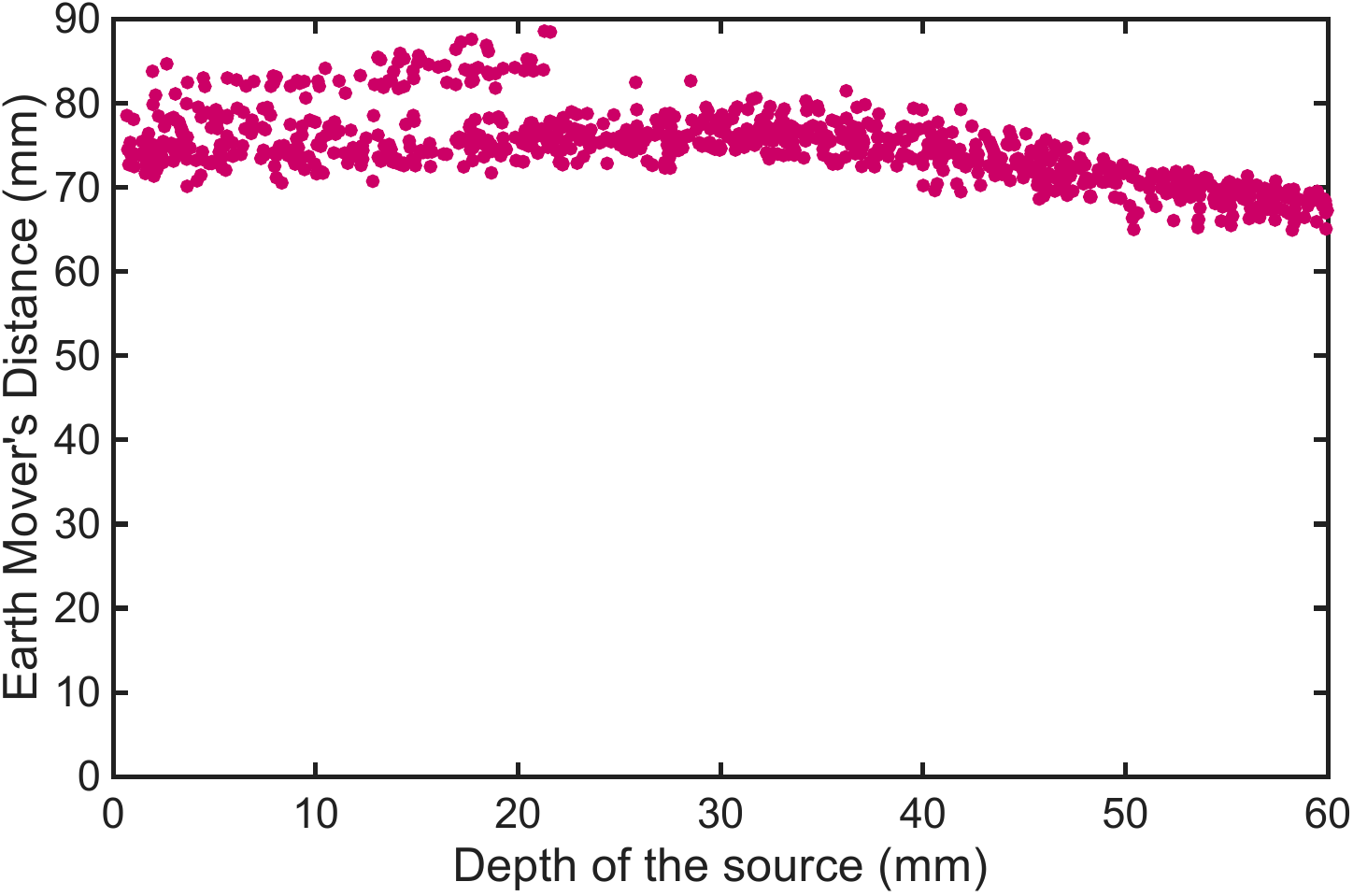}
    \end{minipage}\begin{minipage}{\figWidth}
    \centering
    SHAL1R \vspace{0.1cm}

        \includegraphics[width=0.98\linewidth]{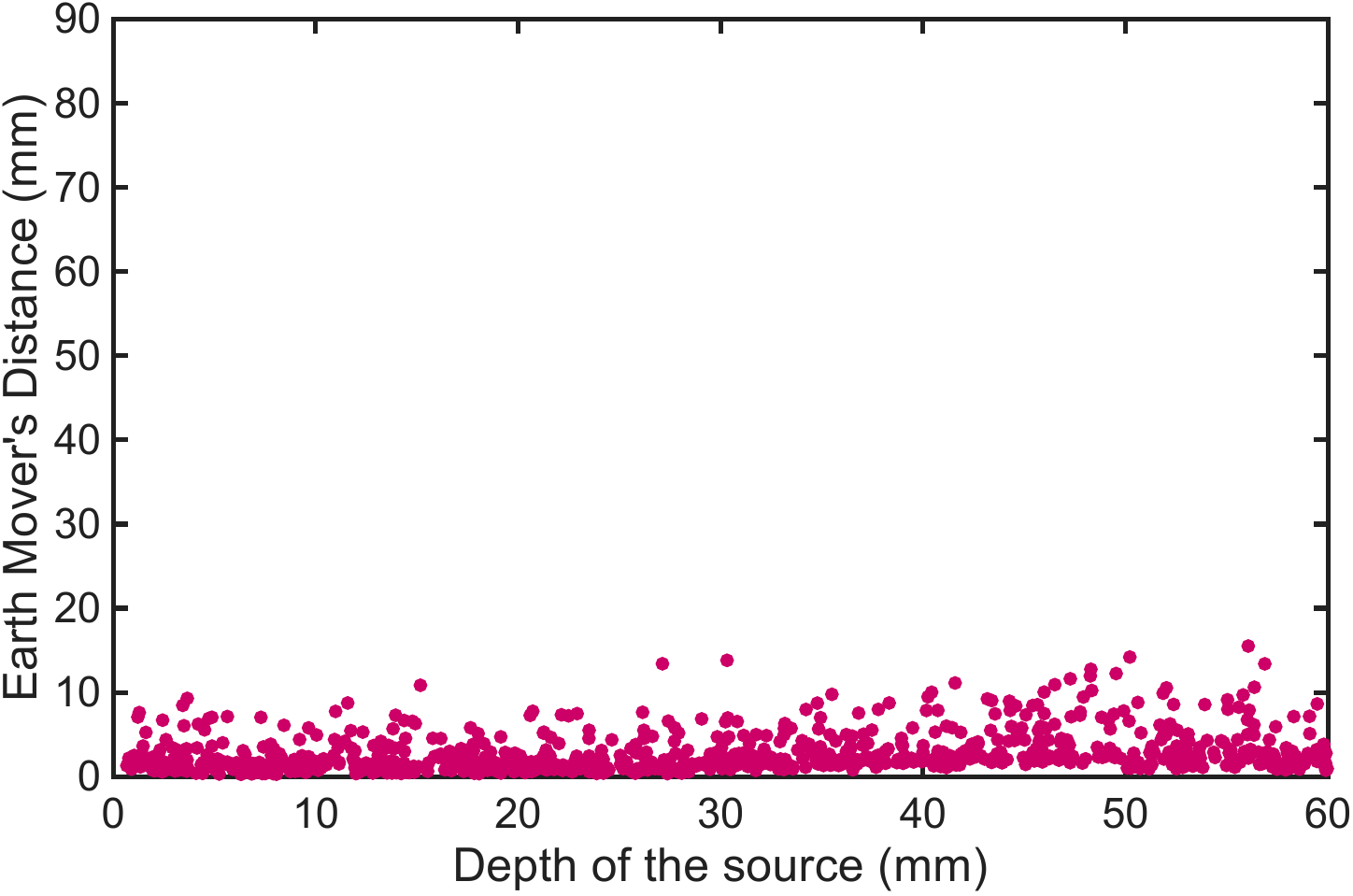}
    \end{minipage}

    \begin{minipage}{0.05\linewidth}
        \rotatebox{90}{$\Hdiv$}
    \end{minipage}\begin{minipage}{\figWidth}
    \centering
        \includegraphics[width=0.98\linewidth]{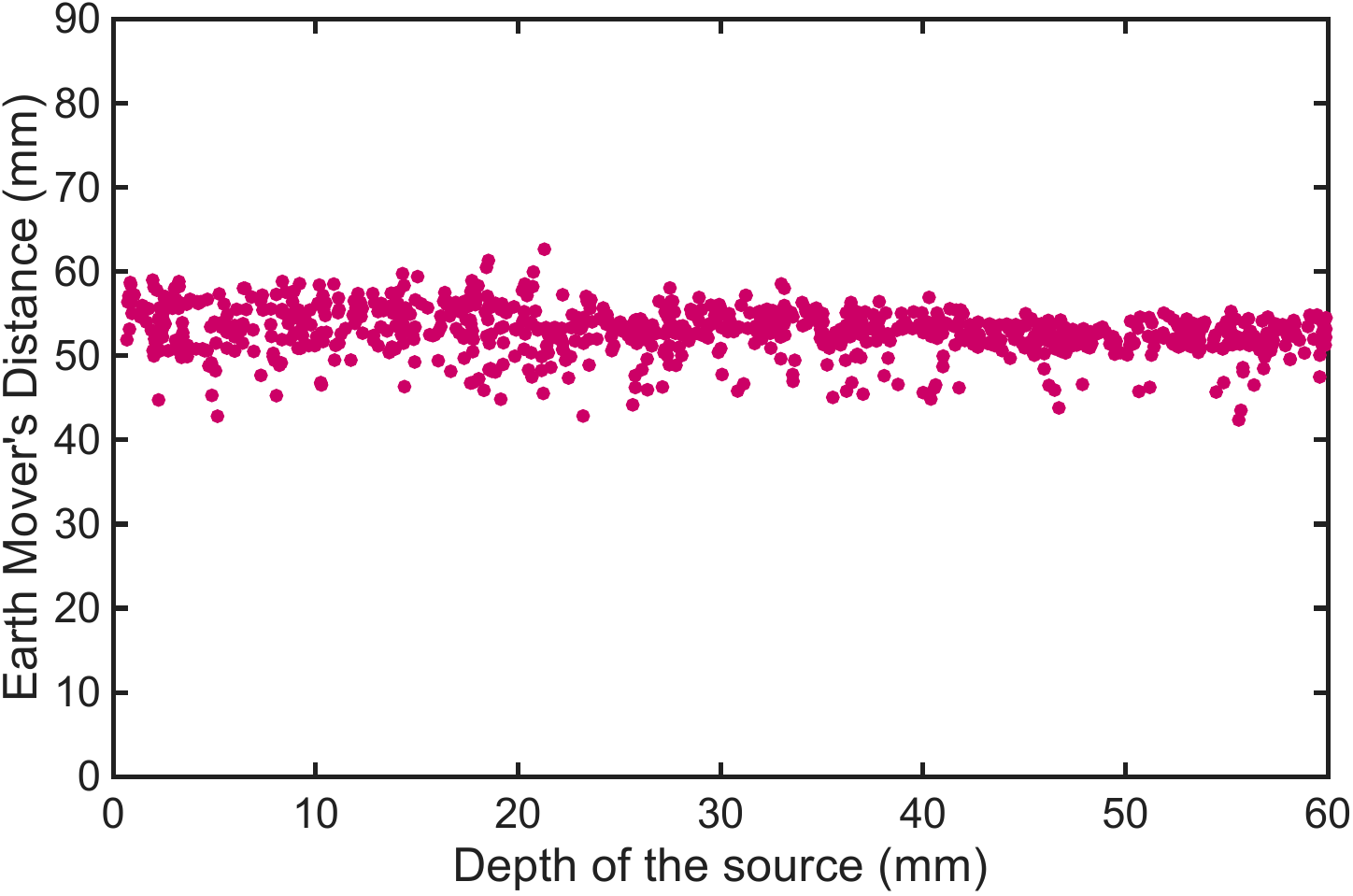}
    \end{minipage}\begin{minipage}{\figWidth}
    \centering
        \includegraphics[width=0.98\linewidth]{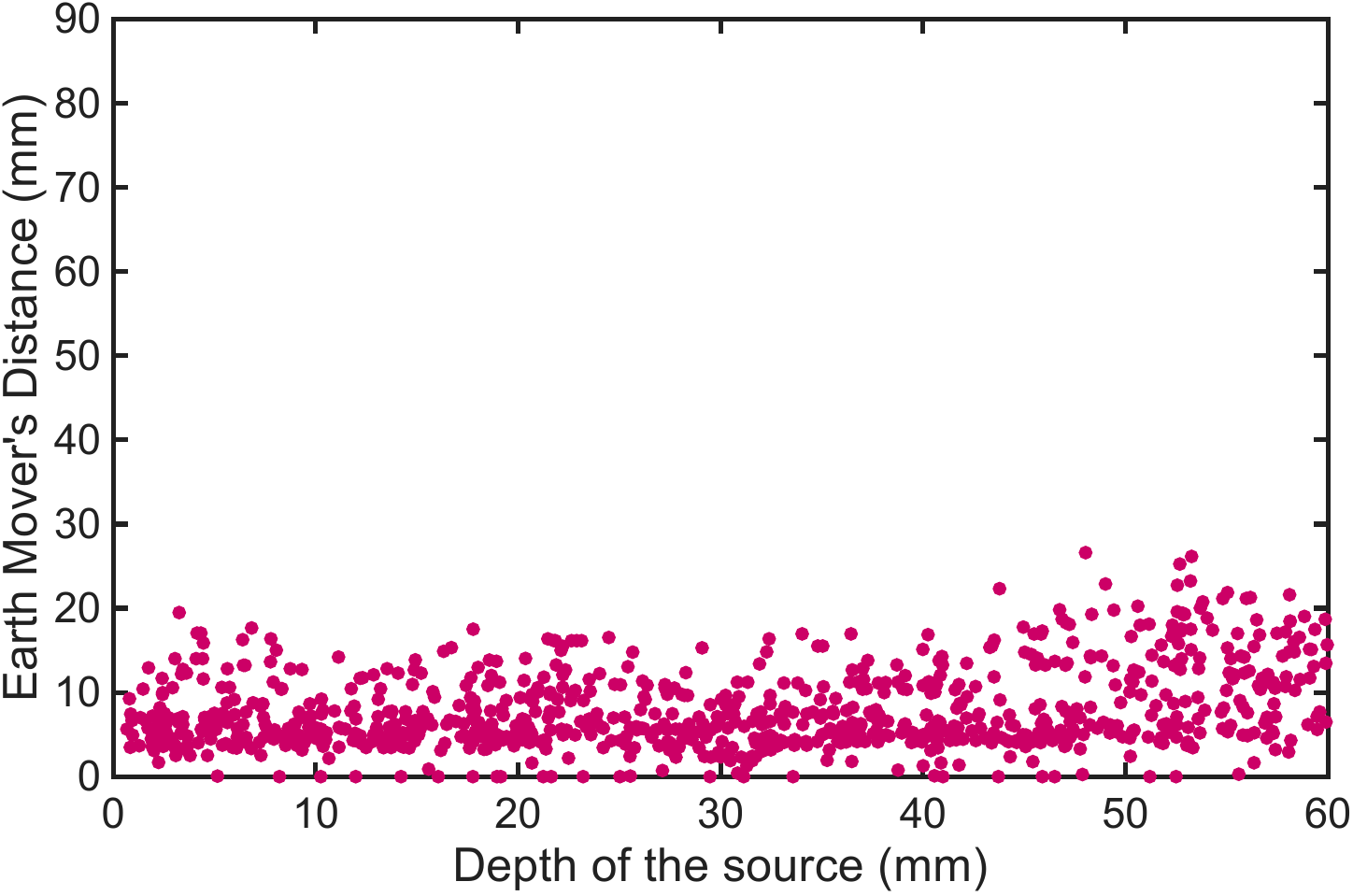}
    \end{minipage}

    \begin{minipage}{0.05\linewidth}
        \rotatebox{90}{Local subt.}
    \end{minipage}\begin{minipage}{\figWidth}
    \centering
        \includegraphics[width=0.98\linewidth]{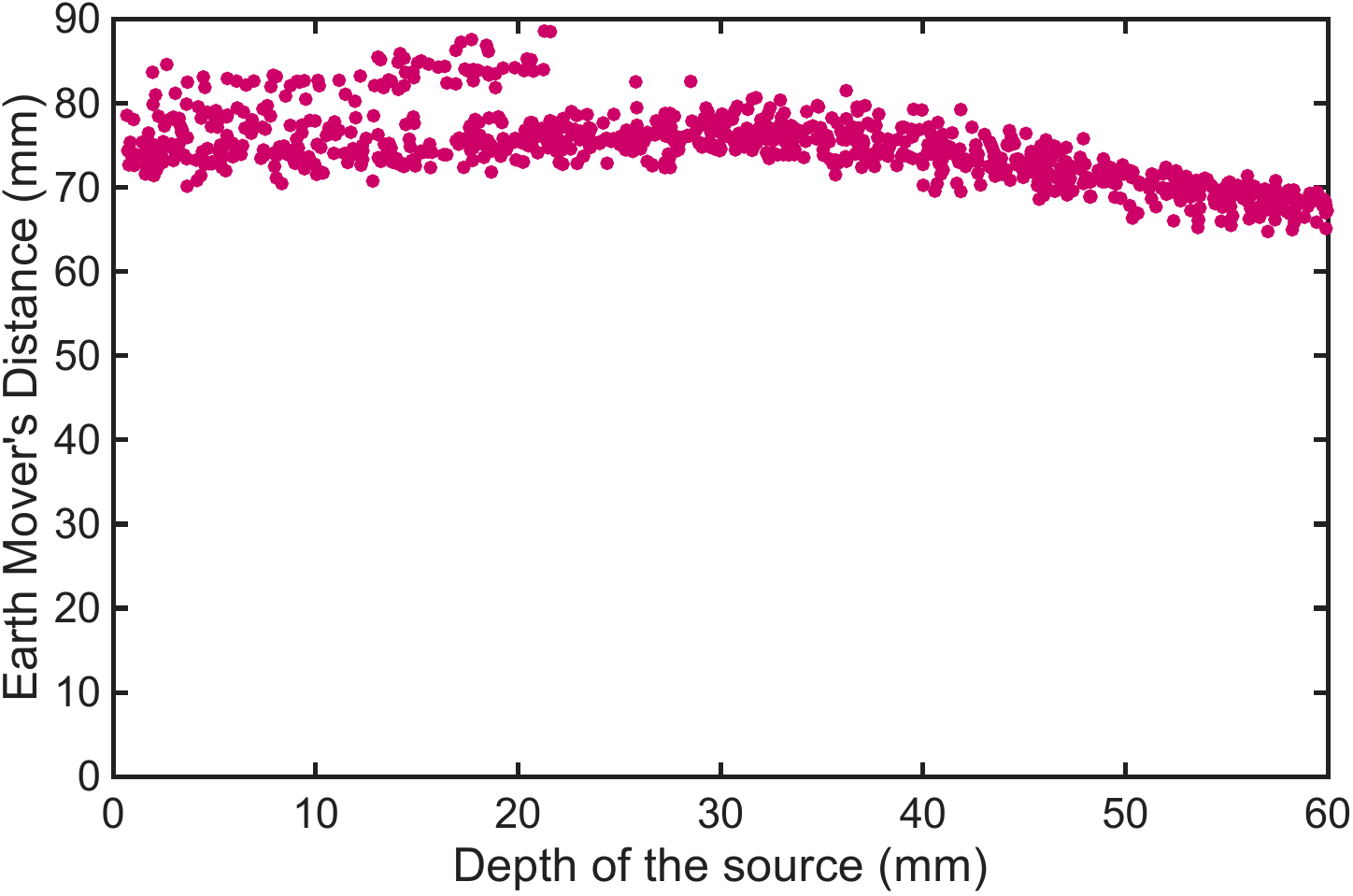}
    \end{minipage}\begin{minipage}{\figWidth}
    \centering
        \includegraphics[width=0.98\linewidth]{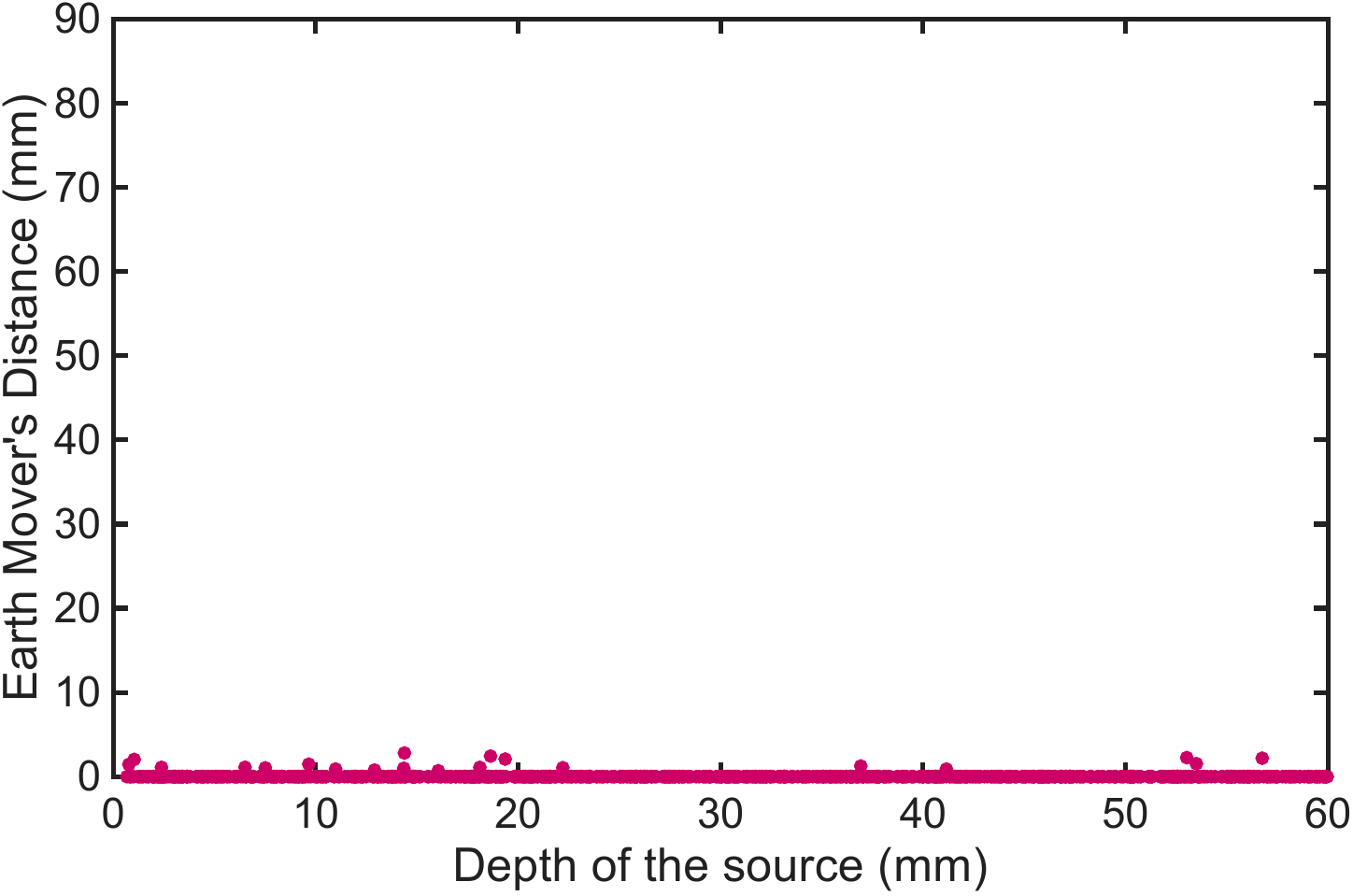}
    \end{minipage}
    \caption{\small Earth Mover's Distances for estimated sources.}
    \label{fig.EMD}
\end{figure}

For a more concrete display of how well each forward and inverse method combination presented in Figure~\ref{fig.depthbias} performs, Figure~\ref{fig.EMD} shows a set of Earth Mover's Distances (EMD) between true and estimated sources at different distances from the inner skull surface.

Here the observation is similar to what it was before: SHAL1R produces a systematically better estimate than sLORETA, regardless of the chosen forward interpolation method, with the bulk of the Whitney and SHAL1R EMDs being near the \qty{2}{\milli\meter} mark, while with Whitney and sLORETA the distances are mostly spread between approximately \num{60}--\qty{80}{\milli\meter}. In the case of $\Hdiv$, the average EMD with SHAL1R is around \qty{5}{\milli\meter}, while majority of the sLORETA's EMD values are between \num{50}--\qty{60}{\milli\meter}. With Local subtraction and sLORETA, we see no significant change in the EMD distribution when compared to Whitney and sLORETA, but a change to Local subtraction clearly improves the correspondence of true and reconstructed depths, as the EMD drops mostly to \qty{0}{\milli\meter} with this change of the forward interpolation method, with a few non-zero outliers present especially with superficial sources. As depth increases, we observe a decreasing trend in EMD with sLORETA.

\section{Discussion}\label{sec.discussion}

We studied EEG source imaging as a coupled forward--inverse modelling problem and used distributional signatures as a way to analyze reconstruction results. The main focus was on studying the effects of alternative source models and their implementations: DUNEuro's \cite{schrader2021duneuro} face-intersecting and edgewise Whitney basis functions, $\Hdiv$ implemented in Zeffiro Interface \cite{he2020zeffiro}, and Local subtraction implemented in DUNEuro. The findings show that different modelling choices lead to distinct activity patterns, even when localization accuracy is similar. Namely, we examined the spread and smoothness of cortical source reconstructions with sLORETA, SHAL1R, SKF, and DS; measured the depth-dependency of standardized methodologies, sLORETA, and SHAL1R; and computed Earth Mover's Distances (EMD) for reconstructions from these two methods.

The choice of Zeffiro Interface \cite{he2020zeffiro} and DUNEuro \cite{schrader2021duneuro} as the forward modelling software was the result of them being one of the few FEM-based software developed specifically with electromagnetic field modelling of the human brain in mind. Alternative interfaces such as Brainstorm \cite{tadel2011brainstorm} and OpenMEEG \cite{gramfort2010openmeeg} exist, but these are either interfaces that invoke other software (including Zeffiro and DUNEuro), to achieve similar results, or which rely on BEMs, which do not take volumetric anatomical features of the domain into account. Based on the forward solution results, DUNEuro seems to perform better in producing a lead field column norm that one would expect: strong only near the electrodes due to the application of reciprocity \cite{rush1969reciprocity,vanrumste2004reciprocity}, and decreasing smoothly as one moves away from the electrodes. In the case of Zeffiro, we see a field that is less localised near the electrodes and which also extends further into the domain even when outliers above the \qty{95}{\percent} quantile in the field strength have been filtered out. The inverse methods we utilized were also implemented as a part of Zeffiro Interface under the MATLAB namespace \texttt{+inverse}.

Based on the reconstructions of the cortical source and depth bias plot, Whitney basis and Local subtraction seem fairly similar. However, EMDs for SHAL1R show a significant difference between these two source models. SHAL1R is an alternative to sLORETA that explicitly assumes a sparse activity distribution \cite{lahtinen2024shalpr}, which is appropriate since the activity originates from a point source. The SHAL1R's EMD, which is almost zero at every depth with Local subtraction, indicates that there is no other dipole that could produce similar data. Reflecting on studies on reconstructing deep brain activity from human EEG recordings, EEG's capability to distinguish deep sources is evident, but the main challenge is in recovering those weak signals from noisy data \cite{KrishnaswamyPavitra2017Seeo,Fahimi2020ECoGvdEEGdeeppossible,Rezaei2021}. Therefore, the ambiguity should come from the modelled noise, not from the lead field in a numerical setting.

In contrast to the results obtained with Whitney basis and Local subtraction, the worst EMD for SHAL1R is obtained with $\Hdiv$. The cortical reconstructions with SHAL1R and DS indicate that this source model yields a region of nearly identical sources at the cortical level; hence, the difficulty of localizing the source is greater than with the Whitney basis and Local subtraction. High EMD values for both sLORETA and SHAL1R indicate high ambiguity in the source location at every depth.

Considering the decreasing tendency of EMD with sLORETA for Whitney and Local subtraction, this indicates that the method has an inflated ability to reconstruct sources farther from the sensors than near them. This could indicate an unrealistic model of deep activity, especially patch-like or spread-out activity, since standardized methods offer higher estimation accuracy for cortical sources, besides their unbiasedness \cite{lahtinen-standardization-2024,giri2025localization}. Moreover, SKF provides a more concise estimate of the cortical activity for $\Hdiv$ than others. On one hand, this further consolidates the notion that a source model with weaker source separation is better suited to inversion models that assume wider distributions, such as sLORETA and SKF. On the other hand, a sparsity-promoting method, such as SHAL1R and DS estimation with the best goodness-of-fit, could be preferred choices when the source to be recovered is focal, i.e., concentrated in a small region, and the source model used in lead field generation supports the notion of point sources.

The study is limited to estimating a single dipole. This is partly because the goal of the paper is to map the differences between forward modelling, and partly due to yet-unknown limitations of the methods for estimating different patch sources and multiple-source configurations. Patch source estimation is one of the future research directions we plan to pursue.

Overall, due to the high performance of SHAL1R and DS, we conclude that Local subtraction yields the most accurate point-source model; however, a source model designed for patch-like activity needs to be developed to reconstruct and estimate sources with spread, as Local subtraction causes sub-optimal behavior for sLORETA.

\section{Author contributions}\label{sec.author.contributions}

Santtu Söderholm produced the EEG lead fields $\leadFieldMatrix$ using DUNEuro and Zeffiro Interface and was one of the main writers of the work. Joonas Lahtinen produced inverse reconstructions of brain activity using the various inverse methods mentioned in this work and also wrote and improved many parts of the text. Sampsa Pursiainen provided the Zeffiro $\Hdiv$ lead field routines, participated in writing the introduction and abstract, and worked in a supervisory role throughout the writing process.

\section{Copyright notice}\label{sec.copyright}

This work has been submitted to the conference MetroXRAINE 2026 as \cite{soderholm2026forwardinverse}. Copyright © 2026 IEEE.


\interlinepenalty=10000

\balance

\printbibliography

\end{document}